\documentclass[12pt]{article}
\usepackage{amssymb}
\usepackage{latexsym}
\newcounter{mytab}

\newtheorem{THM}{Theorem}[section]

\newtheorem{LEM}{Lemma}[section]
\newtheorem{COR}{Corollary}[section]
\begin{document}
%
%

\title{Reflection Groups and Polytopes over Finite Fields, II}
 
\author{B. Monson\thanks{Supported by NSERC of Canada Grant \# 4818}\\
University of New Brunswick\\
Fredericton, New Brunswick, Canada E3B 5A3
\and and\\[.05in]
Egon Schulte\thanks{Supported by NSA-grants H98230-04-1-0116 and
H98230-05-1-0027}\\
Northeastern University\\
Boston, Massachussetts,  USA, 02115}

\date{ \today }
\maketitle

\begin{abstract}
\noindent
When the standard representation of a crystallographic Coxeter group $\Gamma$ is 
reduced modulo an odd prime $p$, a finite representation in some orthogonal space 
over $\mathbb{Z}_p$ is obtained. If $\Gamma$ has a string diagram, the latter group will often be the automorphism group of a finite regular polytope.  In Part I we described the basics of this construction and enumerated the polytopes associated with the groups of rank $3$ and the groups of spherical or Euclidean type. In this paper, we investigate such families of polytopes for more general choices of $\Gamma$, including all groups of rank $4$.
In particular, we study in depth the interplay between their geometric properties and the algebraic structure of the corresponding finite orthogonal group.



\end{abstract}

\section{Introduction}
\label{intro}

The regular polytopes are a rich and ongoing source of 
mathematical ideas.  Their combinatorial features, for instance,
have been beautifully generalized in the theory of 
\textit{abstract regular polytopes}.

In \cite{monsch}, the precursor to this paper, we surveyed some
of the essential properties of an abstract regular polytope $\mathcal{P}$,
referring to \cite{arp} for details. Then, reframing the key results in
\cite{zal1}, we outlined an abbreviated classification of 
finite, irreducible groups generated by reflections
in $n$-space $V$, over a field of odd
characteristic $p$ (see \cite[Thm. 3.1]{monsch}). 

When $G$ is a (possibly infinite) crystallographic Coxeter group with
string diagram, reduction modulo an odd prime $p$ of the standard real representation
yields a finite reflection group $G^p$, which we could then classify 
and which is often the automorphism group of a finite, 
abstract regular $n$-polytope $\mathcal{P}$. (If this is so, we say that $G^p$
is a \textit{string C-group}.)  

Next we established  two useful criteria for $G^p$ to be a string C-group:
Theorems 4.1 and 4.2 of \cite{monsch} concern the features of $V$
as an orthogonal space, as well as the action of standard
subgroups of $G^p$ on $V$. With this, we were able to  classify 
all groups  $G^p$, and their polytopes, whenever $n \leq 3$, as
well as  when $G$ is of spherical or Euclidean type, for
all ranks $n$. 

Here, we begin by summarizing in Section~\ref{basicreview} some key notation. Next, in Section~\ref{interprop}, we extend our  criteria for $G^p$ to be a string C-group. Finally, in Sections~\ref{polfourone} and \ref{polfourtwo},  we discuss and completely classify all  $4$-polytopes which arise from our construction.

\section{Notation}
\label{basicreview}

We refer the reader to the notation and basic set up in \cite{monsch}. Throughout, $G = \langle r_{0},\ldots,r_{n-1} \rangle$ will be a crystallographic Coxeter group $[p_{1},p_{2},\ldots,p_{n-1}]$ with a string Coxeter diagram $\Delta_{c}(G)$ (with branches labeled $p_{1},p_{2},\ldots,p_{n-1}$, respectively), obtained from the corresponding abstract Coxeter group $\Gamma = \langle \rho_{0},\ldots,\rho_{n-1} \rangle$ via the standard representation on real $n$-space~$V$. For any odd prime $p$, we may reduce $G$ modulo $p$ to obtain a subgroup $G^p$ of $GL_n(\mathbb{Z}_p)$ generated by the modular images of the $r_i$'s. We shall abuse notation by referring to the modular images of objects by the same name (such as $r_i$, $b_i$, $B =[b_{ij}]$, $V$, etc.). In particular, $\{b_i\}$ will denote the standard basis for $V = \mathbb{Z}_p^n$.  In any event, $G^p$ is a subgroup of the orthogonal group $O( \mathbb{Z}_p^n)$ of isometries for the (possibly singular) symmetric bilinear form $x \cdot y$,  the latter being defined on $\mathbb{Z}_p^n$ by means of the Gram matrix $B$; in particular, $r_{i}$ is the orthogonal reflection with {\em root\/} $b_i$ if $b_i^2 \neq 0$.  

Next we make a convenient definition: if $p\geq 5$, or $p=3$ but no branch of $\Delta_{c}(G)$ is marked $6$, then we say that $p$ is \textit{generic} for $G$. In such cases, no node label of 
the diagram $\Delta(G)$ (for a basic system) is zero $\bmod\, p$ and a change in the underlying basic system for $G$ has the effect of merely conjugating $G^p$ in $GL_n(\mathbb{Z}_p)$. On the other hand, in the non-generic case, in which $p=3$ and $\Delta_{c}(G)$ has branches marked $6$, the group $G^p$ may depend essentially on the actual diagram $\Delta(G)$ taken for the reduction mod $p$. (Note that $p$ generic does not necessarily mean that $p \nmid |G|$,  or that certain subspaces of $V$ are non-singular, etc.) 

Recall from \cite[Thm. 3.1]{monsch} that an irreducible group $G^p$ of the above sort, 
generated by $n \geq 3$ reflections, must necessarily be one of the following:
\begin{itemize}
\item an orthogonal group $O(n,p,\epsilon) = O(V)$ or $O_j(n,p,\epsilon) = O_{j}(V)$, excluding the cases $O_1(3,3,0)$, $O_2(3,5,0)$, $O_2(5,3,0)$ (supposing for these three that 
$\mbox{disc}(V) \sim 1$), and also excluding the case $O_j(4,3,-1)$; or
\item the reduction mod $p$ of one of the finite linear Coxeter groups of type $A_n$ 
($p\nmid n+1$), $B_n$, $D_n$, $E_6$ ($p \neq 3$), $E_7$, $E_8$, $F_4$, $H_3$ or $H_4$.  
\end{itemize}
We shall say in  these two cases that $G^p$ is of \textit{orthogonal} or \textit{spherical type}, respectively, although there is some overlap for small primes. Our description rests on the classification of the finite irreducible reflection groups over any field, obtained in {Zalesski\u\i} \& {Sere\v zkin} \cite{zal1} (see also \cite{kant,wag1,wag2,zal2}).  It is only a slight abuse of notation to let $[p_1, \ldots , p_{n-1}]^p$ denote the modular representation of a 
group $[p_1, \ldots , p_{n-1}]$, so long as $p$ is generic for the group.

The generators $r_i$ of $G^p$ satisfy the Coxeter-type relations inherited from $G$.  Our main problem is to determine when $G^p$ has the intersection property (\ref{interII}) for its standard subgroups.   For any $J \subseteq \{0, \ldots, n-1\}$, we let 
$G_J^p := \langle r_j \,|\, j \not\in J \rangle $;  in particular, for $k,l \in \{0, \ldots n-1\}$ we let $G_k^p := \langle r_j \,|\, j \neq k\rangle$ and $G_{k,l}^p := \langle r_j \,|\, j \neq k,l \rangle$. 
Then $G^p$ is a \textit{string C-group} if and only if $G^p$ satisfies the 
\textit{intersection property}
\begin{equation}
\label{interII}
\langle \rho_i\,|\, i \in I \rangle \cap \langle \rho_i\,|\,i \in J \rangle =
\langle \rho_i\,|\, i \in I \cap J \rangle ;
\end{equation}
and in this case $G^p$ is the automorphism group of a finite regular polytope 
denoted by ${\cal P}(G^p)$ (see \cite[\S 2E]{arp}).  Note as well that $G^p$ is a string C-group if and only if $G_{0}^p$ and $G_{n-1}^p$ are string C-groups and 
$G_{0}^p \cap G_{n-1}^p = G_{0,n-1}^p$. We also let $V_J$ be the subspace of 
$V = \mathbb{Z}_p^n$ spanned by $\{b_j \!\mid\! j \not\in J \}$, and similarly for $V_k, V_{k,l}$. Note that $V_J$ is $G_J^p$-invariant.

\section{The Intersection Property}
\label{interprop}

The goal of this section is to assess when the modular reduction $G^p$ of a crystallographic 
string Coxeter group $G$ satisfies the intersection property (\ref{interII}). In \cite{monsch} we 
established a number of sufficient conditions and verified that $G^p$, with $p\geq 3$, has the 
intersection property whenever $G$ has rank at most $3$, or whenever $G$ is of spherical or Euclidean type. However, the situation changes drastically for more general groups of higher ranks, with obstructions already occurring for rank $4$.  We already know that $G^p$ is a string C-group if one of the subgroups $G_{0}^p$ or $G_{n-1}^p$ is spherical and the other is a string C-group (see \cite[Thm. 4.2]{monsch}). Moreover, $G^p$ also is a string C-group if both $G_{0}^p$ and $G_{n-1}^p$ are C-groups, $V_{0,n-1}$ is a non-singular subspace of $V$, and $G_{0,n-1}^p$ is the full orthogonal group $O(n-2,p,\epsilon)$ on $V_{0,n-1}$ (see \cite[Thm. 4.1]{monsch}). 

The criteria established here will settle the Coxeter groups $[k,l,m]$ of rank $4$ completely. Before we move on, note two simple cases. If $l=2$ and $k,m < \infty$, then 
\[ G \cong [k] \times [m] \cong G^{p} , \]
so $G^p$ certainly is a C-group. Similarly, if say  $m = 2$ (but $k$ or $l = \infty$ is allowed), then $G \simeq [k,l] \times C_2$, and
$$ G^p \simeq [k,l]^p \times C_2\;.$$
Thus the intersection property  of $G^p$ follows directly from that of its subgroup $[k,l]^p$ 
(see \cite[Thm. 5.1]{monsch}). More generally, if $[p_1, \ldots , p_{n-1}]^p$ is a string C-group, then so is
$$[p_1, \ldots , p_{n-1}, 2]^p \simeq  [p_1, \ldots , p_{n-1}]^p \times C_2\;\;. $$

After a preliminary lemma, we shall continue to build upon the known results
concerning $C$-groups mentioned above.

\noindent
\begin{LEM}\label{3new1} 
Suppose that $G$ has rank $n$ and that  
$ (\mbox{rad}\;V) \cap V_j = \{o\}$. 
Then the subgroup $G_j^p$ is, 
by restriction to the  invariant subspace $V_j$, 
isomorphic to $H^p$,  the reduction modulo $p$ of the group of rank $n-1$ 
defined  from the subdiagram of $\Delta(G)$ which results from the 
deletion of node $j$. In particular, when $j=0$ or $n-1$
this holds if  $p$ is generic for  
 $G$.  
\end{LEM}

\noindent
\textbf{Proof}. This result is well known in characteristic $0$
\cite[\S 5.5]{humph}. Here we restrict $g \in G_j^p$ to the invariant subspace $V_j$,
and so obtain a homomorphism
\begin{eqnarray*}
\varphi: & G_j^p & \longrightarrow O(V_j)\\
         & g   & \longmapsto g_{|V_{j}}\;\;.
\end{eqnarray*}
Of course, as a subspace of $V$, $V_j$  is isometric 
to  $\mathbb{Z}_p^{n-1}$, with the metric structure obtained
from the subdiagram   of $\Delta(G)$  obtained by deleting   node $j$.
Clearly the image group 
$\varphi(G_j^p)$  is isomorphic to the reflection group $H^p$ of rank $n-1$
defined directly from the subdiagram.

Suppose $g \in \ker \varphi$.  Then $g(b_k) =
b_k$ for all $k \neq j$, whereas $g(b_j) = b_j + x$ for some $x \in
V_j$.  Thus for any $k \neq j$
$$
b_j \cdot b_k  =  g(b_j) \cdot g(b_k)
	=  (b_j + x) \cdot b_k =  b_j \cdot b_k +  x \cdot b_k\;,
$$
so that $x \cdot b_k = 0$, and $x \in \mbox{rad}\;V_j$.  But then $x
\cdot x = 0$ and so
$$ b_j \cdot b_j  =  g(b_j) \cdot g(b_j)
	=  b_j \cdot b_j + 2x \cdot b_j + x \cdot x\;,
$$
whence $x \cdot b_j = 0$.  Thus $x \in (\mbox{rad}\;V) \cap V_j$, so that $x=o$  when this 
subspace is trivial. Hence $\varphi$ is injective. When $p$ is generic for $G$, a direct calculation in coordinates along the string diagram shows that   $(\mbox{rad}\;V) \cap V_j 
= \{o\}$ for $j=0, n-1$.
\hfill$\square$

\noindent
\textbf{Remark}. Informally,  the Lemma asserts that reduction by a generic prime
commutes with the deletion of a node from $\Delta(G)$. Note that
$$G_j^p \simeq [p_1, \ldots , p_{j-1}]^p \times 
  [p_{j+2}, \ldots , p_{n-1}]^p \simeq 
  [p_1, \ldots , p_{j-1}, 2,p_{j+2}, \ldots , p_{n-1} ]^p\;\;.$$

Concerning the non-generic cases,  there are examples showing the
necessity of the hypotheses. For example, the group 
$G \simeq [4,3,6]$ with diagram 
$$ \stackrel{2}{\bullet}\frac{}{\;\;\;\;\;\;\;}
	\stackrel{1}{\bullet}\frac{}{\;\;\;\;\;\;\;}
	\stackrel{1}{\bullet}\frac{}{\;\;\;\;\;\;\;}
	\stackrel{3}{\bullet} 
$$
yields,  as we observe below, a $C$-group $G^3$. Here
the subgroup $G_0^3$ is the automorphism group 
of order $108$ for the toroidal polyhedron $\{3, 6\}_{(3,0)}$. However,
the subdiagram
$$\stackrel{1}{\bullet}\frac{}{\;\;\;\;\;\;\;}
	\stackrel{1}{\bullet}\frac{}{\;\;\;\;\;\;\;}
	\stackrel{3}{\bullet} $$
yields the smaller  group of order $36$ for $\{3, 6\}_{(1,1)}$. 
Thus the map $\varphi$ 
of the Lemma is here not injective.

\begin{THM}
\label{thmfourthree}
Suppose that $G \simeq [k,l,m]$ is crystallographic
and that the subgroup $[k,l]$ or $[l,m]$ is spherical. Then $G^p$ is a C-group 
for any prime $p \geq 3$.
\end{THM}

\noindent
\textbf{Proof.}  Let $[k,l]$ (say) be spherical, so that $G_3^p \simeq G_3 = [k,l]$. 
First suppose that $p$ is generic for $G$. Since $G_0^p$ is a C-group by Lemma~\ref{3new1}, the proof follows directly from \cite[Thm. 4.2]{monsch}.
Moreover, even in non-generic cases of rank $4$ (so that $p=3$), 
$G^3$ turns out to be  a C-group when $[k,l]$ is spherical. This is routinely verified using
the computer algebra system GAP \cite{gap}. The 
pertinent examples are $G \simeq [3,3,6], [3,4,6]$ or $[4,3,6]$, each with two essentially distinct diagrams $\Delta(G)$ (for the basic systems). 
\hfill $\square$

We now establish two general results, of which the first allows us to reject large classes of 
groups $G^p$ as C-groups because of the size of their subgroups $G_{0}^p \cap G_{n-1}^p$.  
First we deal with the \textit{fully non-singular} case.

\begin{THM}
\label{midnondeg} 
Let $G = \langle r_0, \ldots , r_{n-1} \rangle$  be a crystallographic linear Coxeter group with 
string diagram. Suppose  that $n\geq 3$ and  that the prime $p$ is generic for $G$.
Let the subspaces $V$, $V_{0}$, $V_{n-1}$ and  
$V_{0,n-1}$ be non-singular, and let $G_0^p$, $G_{n-1}^p$ be of orthogonal type. 
Suppose as well that there is a square among the labels of the nodes $1,\ldots,n-2$ of the 
diagram $\Delta (G)$ (this can be achieved by readjusting the node labels). \\ 
(a) Then $G_{0}^p \cap G_{n-1}^p$ acts trivially on $V_{0,n-1}^{\perp}$, and
\[ O_{1}(V_{0,n-1}) \leq G_{0}^p \cap G_{n-1}^p \leq  O (V_{0,n-1}) , \]
where we have identified the $(n-2)$-dimensional groups $O (V_{0,n-1})$ and $O_1(V_{0,n-1})$ with the pointwise stabilizers of $V_{0,n-1}^{\perp}$ in the $n$-dimensional groups $O (V)$ and $O_1(V)$, respectively.\\[.01in]
(b) If $G_{0}^{p} = O(V_{0})$ and $G_{n-1}^{p} = O(V_{n-1})$, with a similar interpretation as 
stabilizers, then 
\[ G_{0}^p \cap G_{n-1}^p = O (V_{0,n-1}) . \]
\end{THM}

\noindent
\textbf{Proof.}  
Since all four subspaces of $V = \langle b_{0},\ldots,b_{n-1} \rangle$ are non-singular, we have the orthogonal sums 
\[ V = V_{n-1} \oplus \langle v \rangle = V_{0} \oplus \langle v' \rangle, \;
V_{n-1} = V_{0,n-1} \oplus \langle w \rangle, \;
V_{0} = V_{0,n-1} \oplus \langle w' \rangle, \]
for non-isotropic vectors $v,v',w,w'$. Then,
\[ \langle v,v' \rangle = V_{0,n-1}^{\perp} = \langle w,w' \rangle . \]
Since $p$ is generic for $G$, each reflection $r_j$ actually has $b_j$ as a root, and
$v \perp b_{j}$ for $j \leq n-2$, while $v' \perp b_{j}$ for 
$j \geq 1$.  Hence the subgroups $G_{n-1}^p$, $G_{0}^p$ and $G_{0}^p \cap G_{n-1}^p$ stabilize the vectors $v$, $v'$ or $v,v'$, respectively.  In particular,  
\[ G_{0}^p \cap G_{n-1}^p \leq O(V_{0,n-1}) , \]
with $O(V_{0,n-1})$ identified with the pointwise stabilizer of $V_{0,n-1}^{\perp}$ in 
$O (V)$. Note that the restrictions of $G_{n-1}^p$, $G_{0}^p$ and $G_{0}^p \cap G_{n-1}^p$ to the subspaces $V_{n-1}$, $V_{0}$ or $V_{0,n-1}$, respectively, are faithful, by Lemma~\ref{3new1}.

Since $G_0^p$, $G_{n-1}^p$ are  of orthogonal type 
and there is a square among the labels of the 
nodes $1,\ldots,n-2$,  we must have $O_{1}(V_{n-1}) \leq G_{n-1}^p$ and 
$O_{1}(V_{0}) \leq G_{0}^p$ (that is, a group merely of type $O_{2}$ cannot occur). Now, if 
$g \in O_{1}(V_{0,n-1})$, then $g(v)=v$, so that $g  \in O(V_{n-1})$; but the spinor norm is 
invariant under orthogonal embedding (\cite[Thm. 5.13]{art}), so actually 
$g  \in O_{1}(V_{n-1})$. Similary, $g  \in O_{1}(V_{0})$, and hence 
$g \in G_{0}^p \cap G_{n-1}^p$.  This completes the proof of  part (a).

Now let $G_{0}^{p} = O(V_{0})$ and $G_{n-1}^{p} = O(V_{n-1})$. Once again, if 
$g \in O(V_{0,n-1})$, then $g(v)=v$, so now $g  \in O(V_{n-1}) = G_{n-1}^{p}$. Similarly, 
$g  \in O(V_{0}) = G_{0}^{p}$, and hence $g \in G_{0}^p \cap G_{n-1}^p$, as required.
\hfill $\square$

We note an immediate corollary to Theorem~\ref{midnondeg}. It shows that many groups 
$G^p$ of rank $4$ fail to satisfy the intersection property for large primes $p$. However, 
those primes for which $G^p$ actually is a C-group lead to interesting polytopes, which we  
investigate in later sections. 

\begin{COR}
\label{corfourone}
Suppose the prime $p$ is generic for the crystallographic group  $G =[k,l,m]$. 
Let  $V$, $V_{0}$, $V_{3}$, $V_{0,3}$ be non-singular, and let 
$G_{0}^p$, $G_{3}^p$ be of orthogonal type. \\ 
(a) Then $G^p$ is not a C-group if $p > 2l + \epsilon(V_{0,3})$, where 
$\epsilon(V_{0,3}) = \pm 1$ is the parameter associated with the plane $V_{0,3}$.\\
(b) If $G_{0}^{p} = O(V_{0})$ and $G_{3}^{p} = O(V_{3})$, then $G^p$ is a C-group if and only 
if $p = l + \epsilon(V_{0,3})$.
\end{COR}

{\bf Proof.}  
We apply Theorem~\ref{midnondeg} with $n=4$.  The subgroups $G_{0}^p$ and $G_{3}^p$ are known to be C-groups (\cite[Thm. 5.1]{monsch}), so it suffices to determine when
\[ G_{0}^p \cap G_{3}^p = G_{0,3}^p . \]   
Now $G_{0,3}^p = \langle r_{1},r_{2} \rangle$ is a dihedral group of order $2l = 6$, $8$ or 
$12$; the case $l = \infty$ is excluded, as $V_{0,3}$ is then a non-singular plane.  Note that we may assume that there is a square (in fact, a $1$) among the labels of the nodes $1$ or $2$ of the diagram; this can be achieved by readjusting the node labels as described earlier.  Then, by Theorem~\ref{midnondeg} we have 
\[O_{1}(V_{0,3}) \leq G_{0}^p \cap G_{3}^p , \] 
so $G_{0}^p \cap G_{3}^p$ is larger than $G_{0,3}^p$ if the order of $O_{1}(V_{0,3})$, which 
is  $p-\epsilon(V_{03})$, exceeds $2l$. Hence the intersection property certainly fails if 
$p > 2l + \epsilon(V_{0,3})$.  Moreover, by Theorem~\ref{midnondeg}, if 
$G_{0}^{p} = O(V_{0})$ and $G_{3}^{p} = O(V_{3})$, then 
\[ G_{0}^p \cap G_{3}^p = O(V_{0,3}) , \] 
so $G_{0}^p \cap G_{3}^p = G_{0,3}^p$ if and only if $2(p-\epsilon(V_{03})) = 2l$, or 
equivalently, $p = l + \epsilon(V_{0,3})$.
\hfill $\square$

Corollary~\ref{corfourone} immediately implies (in fully non-singular cases)
that $G^p$ is not a C-group if $p > 13$. 
However, for the primes $p = 5,7,11,13$ (and $3$), the outcome is less predictable and 
actually depends on the group $G = [k,l,m]$ as well as the diagram $\Delta(G)$ chosen for 
the reduction modulo $p$. For example, $G^{13}$ can only be a C-group if $l=6$ and 
$G_{0}^{13} = O_1(V_{0})$ or $G_{3}^{13} = O_1(V_{3})$.  Similarly, if $G_{0}^{p} = O(V_{0})$ and $G_{3}^{p} = O(V_{3})$, then $G^p$ is not a C-group if $p > 7$; moreover, $G^7$ can then be a C-group only if $l=6$.  

Next we study the case when the middle section of the diagram for $G$ determines a singular 
space $V_{0,n-1}$, while $V$, $V_{0}$ and $V_{n-1}$ still are non-singular, again with 
$G_0^p, G_{n-1}^p$ of orthogonal type. 
In a singular space $W$ over a field $\mathbb K$, the  isometry  group $O(W)$ leaves invariant the radical subspace $\mbox{\rm rad}\, W$, thereby providing a natural 
epimorphism $\eta: O(W) \rightarrow O(W/\mbox{\rm rad}\, W)$. Since
$W/\mbox{\rm rad}\, W$ is non-singular, we may define a `spinor norm' $\theta$
on $W$, sufficient for our needs, by
$$\theta(g):= \theta_{W/\mbox{\rm rad}\, W}(\eta(g)),\;\;g \in O(W)\;. $$ 
Now let $\widehat{O}(W)$ denote the subgroup of $O(W)$ 
consisting of those isometries $g$ which act trivially on $\mbox{\rm rad}\, W$. 
It is not hard to show that $\widehat{O}(W)$  contains  all reflections 
(with non-isotropic roots in $W$) and is even generated by them.
The key observation here is that any transvection in the kernel of
the action of $O(W)$ on $\mbox{\rm rad}\, W$ can be factored as 
a product of reflections (cf. \cite[Thm. 3.20 and p. 133]{art}). 

It is easy to see that if $g$ is the product of reflections with non-isotropic roots 
$a_1, \ldots , a_k$, then $\theta(g) = a_1^2 \cdots a_k^2\, \dot{\mathbb{K}}^2$.
Naturally, by $\widehat{O}_{1}(W)$ (or $\widehat{O}_{2}(W)$) we mean the subgroup of $\widehat{O}(W)$ generated by the 
reflections in $O(W)$ whose spinor norm is a square (or non-square, respectively).  

\begin{THM}
\label{middeg} 
Let $G = \langle r_0, \ldots , r_{n-1} \rangle$  be a crystallographic linear Coxeter group 
with string diagram, and suppose the prime $p$ is generic for $G$. Let $V$, $V_{0}$, $V_{n-1}$ be non-singular, let $V_{0,n-1}$ be singular, and let $G_0^p, G_{n-1}^p$ be of orthogonal type. Suppose there is a square among the labels of the nodes $1,\ldots,n-2$ of the diagram $\Delta(G)$ (this can be achieved by readjusting the node labels). \\ 
(a) Then $G_{0}^p \cap G_{n-1}^p$ acts trivially on $V_{0,n-1}^{\perp}$, and
\[ \widehat{O}_{1}(V_{0,n-1}) \leq G_{0}^p \cap G_{n-1}^p \leq  \widehat{O} (V_{0,n-1}) , \]
where $\widehat{O} (V_{0,n-1})$ has been identified with the pointwise stabilizer of 
$V_{0,n-1}^{\perp}$ in $O (V)$, and $\widehat{O}_{1} (V_{0,n-1})$ with the subgroup of 
$O_{1}(V)$ generated
by the reflections with roots in $V_{0,n-1}$ and square spinor norm. \\ 
(b) If $G_{0}^{p} = O(V_{0})$ and $G_{n-1}^{p} = O(V_{n-1})$, then also
\[ \widehat{O}(V_{0,n-1}) = G_{0}^p \cap G_{n-1}^p . \] 
\end{THM}

\noindent
\textbf{Proof.}  
As before we have 
\[V = V_{n-1} \oplus \langle v \rangle = V_{0} \oplus \langle v' \rangle \] 
with non-isotropic vectors $v,v'$.  The subspace $V_{0,n-2}^{\perp}$ still is $2$-dimensional,
 so necessarily $V_{0,n-1}^{\perp} = \langle v,v' \rangle$.  Moreover,  
\[ V_{0,n-1} \cap V_{0,n-1}^{\perp} = \mbox{\rm rad}\, V_{0,n-1}  \neq \{o\} , \]
so the vectors in $V_{0,n-1} \cup V_{0,n-1}^{\perp}$ span a singular hyperplane $U$ in $V$ with $1$-dimensional radical ${\rm rad}\,U = {\rm rad}\, V_{0,n-1}$.  For the same reason as before, $G_{n-1}^p$, $G_{0}^p$ and $G_{0}^p \cap G_{n-1}^p$ stabilize the vectors $v$, $v'$ or $v,v'$, respectively, and yield faithful restrictions to the subspaces $V_{n-1}$, $V_{0}$ and $V_{0,n-1}$. In particular,
\[ G_{0}^p \cap G_{n-1}^p \leq H:= 
\{ g \in O(V) \!\mid\! g(x)=x,  \;\forall x \in V_{0,n-1}^{\perp} \}. \] 
Note also that $G_{0}^p \cap G_{n-1}^p$ leaves $U$ invariant because it leaves $V_{0,n-1}$ 
invariant.

We claim that we may identify $H$ with $\widehat{O}(V_{0,n-1})$; this would 
settle the inclusion on the right in part (a) of the theorem.  Now, since each element of $H$ 
leaves $V_{0,n-1}$ invariant, while fixing $\mbox{rad}\, V_{0,n-1}$, 
we can consider the restriction mapping to $V_{0,n-1}$,
\begin{equation}
\label{kap}
\begin{array}{rcl}
\kappa:\ H& \longrightarrow & \widehat{O}(V_{0,n-1})  \\
g   & \longrightarrow & g_{\big | V_{0,n-1}} .
\end{array}
\end{equation}
We prove that $\kappa$ is an isomorphism. If $g\in \ker(\kappa)$, then $g$ 
acts trivially on $V_{0,n-1}$, and hence also on $U$ because $g\in H$.  It follows 
that $g=e$, the identity mapping on $V$. Here we have used the fact that 
an isometry of a non-singular space is uniquely determined by 
its effect on a hyperplane, $H$ in this case, if this hyperplane is singular 
(\cite[Thm. 3.17]{art}). 
This shows that $\kappa$ is injective. Now let $h \in \widehat{O}(V_{0,n-1})$. 
Then $h$ acts trivially on 
$\mbox{\rm rad}\, V_{0,n-1}$, 
so we can extend $h$ to an isometry $h'$ (say) of 
$U = V_{0,n-1} \oplus \langle v \rangle $ 
by setting $h'(v):=v$. We now apply Witt's extension theorem for 
isometries between 
subspaces of a non-singular space (\cite[Thm. 3.9]{art}) and conclude that 
$h'$ extends further 
to an isometry $g$ of the entire space $V$.  Then $g$ must be in $H$ and 
so $\kappa$ is also 
surjective. By our earlier remarks, $ \widehat{O}(V_{0,n-1})$ is generated 
by all reflections $r_a$, for non-isotropic roots  
$a \in V_{0,n-1}$; pulling back, a similar  claim is true for $H$. 

Continuing along these lines, we now prove  
the inclusion on the left in part (a) of the theorem. For a non-isotropic vector 
$a\in V_{0,n-1}$, let $r_{a,V}$, $r_{a,V_{n-1}}$, $r_{a,V_{0}}$ and $r_{a,V_{0,n-1}}$ denote the reflections with root $a$ in $V$, $V_{n-1}$, $V_{0}$ or $V_{0,n-1}$, respectively.  Then 
$r_{a,V} \in H$ because $a \perp V_{0,n-1}^{\perp}$, and 
\[r_{a,V_{0,n-1}} = \kappa(r_{a,V}).\]  
It follows that the subgroup $H_{1}$ of $H$ generated by the reflections $r_{a,V}$, with 
$a \in V_{0,n-1}$ and $a^{2}$ a square, is isomorphic, under $\kappa$, to 
$\widehat{O}_{1}(V_{0,n-1})$. We have to show that $H_{1} \leq G_{0}^p \cap G_{n-1}^p$.  

Now, by assumption, the subgroups $G_0^p , G_{n-1}^p $ of $G$ are of orthogonal type, and there is a square among the labels of the nodes $1,\ldots,n-2$ of the diagram.  It follows that 
$G_{0}^p$ and $G_{n-1}^p$, when restricted to the subspaces $V_{0}$ or $V_{n-1}$, respectively, must contain the groups $O_{1}(V_{0})$ or $O_{1}(V_{n-1})$.  Hence, if we identify the restricted groups with the stabilizers of $v'$ or $v$, respectively, in $O_{1}(V)$, then we see that $O_{1}(V_{0})$ and $O_{1}(V_{n-1})$ are actually subgroups of $G_{0}^p$ and $G_{n-1}^p$. In particular, if $a \in V_{0,n-1}$ and $a^{2}$ is a square, then $r_{a,V}$ 
belongs to $O_1(V_j)$, for $j = 0, n-1$, so $r_{a,V} \in  G_{0}^p \cap G_{n-1}^p$. Now it follows that $H_{1}$ is a subgroup of $G_{0}^p \cap G_{n-1}^p$.  This settles part (a).

Finally, suppose that $G_{0}^{p} = O(V_{0})$ and $G_{n-1}^{p} = O(V_{n-1})$.  Much of the same analysis carries over, but now it is applied to all reflections, including those whose spinor 
norm is a non-square.  In particular, if $a$ is any non-isotropic vector in $V_{0,n-1}$, then 
$r_{a,V} \in G_{0}^p \cap G_{n-1}^p$.  Hence we also have 
\[ \widehat{O}(V_{0,n-1}) \leq G_{0}^p \cap G_{n-1}^p\;, \]
since $\widehat{O}(V_{0,n-1})$ is generated by its reflections. Now part (a) gives the
equality of these  groups. 
\hfill $\square$

\noindent
\textbf{Remark}.  It is not difficult to show that
$ \widehat{O}_1(V_{0,n-1})$ can also be identified with the 
pointwise stabilizer of $ V_{0,n-1}^\perp$ in $O_1(V)$.

For a crystallographic Coxeter group $[k,l,m]$, the middle section of the diagram determines a 
singular subspace if and only if $l=\infty$. In this case, the reduced group $G^p$ is 
always a $C$-group:

\begin{COR}
\label{corfourtwo} 
Let $G \simeq [k,\infty,m]$ be crystallographic.  Then $G^p$ is a 
C-group for any prime $p \geq 3$.
\end{COR}

\noindent
\textbf{Proof.}  
We know that $G_{0}^p$ and $G_{3}^p$ are C-groups (\cite[Thm. 5.1]{monsch}), so again it suffices 
to check that 
\[ G_{0}^p \cap G_{3}^p = G_{0,3}^p . \]
Suppose for the moment that $p$ is generic for $G$, 
so that we may  apply 
Theorem~\ref{middeg} with $n=4$. Then, with few exceptions, $V$ is still non-singular.
Moreover, $V_{0}$ and $V_{3}$ correspond to the subgroups $[\infty,m]$ and $[k,\infty]$ of 
$G$, and hence are known to be non-singular as well 
(except in the non-generic case with $p=3$,
$k,m=6$; see \cite[Sect. 5]{monsch}). 
However, $V_{0,3}$ is a singular plane, so Theorem~\ref{middeg} implies that
\[ \widehat{O}_{1}(V_{0,3}) \leq G_{0}^p \cap G_{3}^p \leq  \widehat{O} (V_{0,3}) ,\]
provided $V$ is non-singular. If the labels of the nodes $1$ and $2$ of the diagram are $1$ and $4$, respectively, as we may assume, then $V_{0,3}$ is a singular plane in which the squared norm of each non-isotropic vector is a square.  In particular, 
\[ \widehat{O}_{1}(V_{0,3}) = \widehat{O} (V_{0,3}) \cong [p] \cong G_{0,3}^p ,\]
and hence $G_{0}^p \cap G_{3}^p = G_{0,3}^p$.  In fact, since $2b_{1}+b_{2}$ spans 
$\mbox{\rm rad}\, V_{0,3}$, a matrix representing an element of $\widehat{O} (V_{0,3})$ in the 
basis $b_{1},2b_{1}+b_{2}$ must necessarily have the form
\[ \left[ \begin{array}{cc}
\pm 1 & 0\\
\mu & 1
\end{array} \right] 
\quad (\mu \in \mathbb{Z}_p )  ,\]
so there are at most $2p$ of them. On the other hand, the restrictions of $r_{1}$ and $r_{2}$ to $V_{0,3}$ already generate a dihedral group $[p]$ contained in $\widehat{O}_{1}(V_{0,3})$, so the three groups must coincide.

The space $V$ is singular for the following groups $G$ (up to duality) and primes 
$p\,(>3)$:\  $[4,\infty,3]$ and 
$[6,\infty,6]$, with $p=5$; $[3,\infty,3]$ and $[4,\infty,6]$, with $p=7$; and $[3,\infty,6]$, 
with $p=13$. In each of these cases, as well as for $p=3$ for  any of the 
crystallographic groups $[k, \infty, m]$, computations in GAP confirm that $G^p$ also is a 
C-group. 
\hfill $\square$

We now concentrate entirely on the groups $G=[k,l,m]$ which are not yet covered by the 
previous results. These groups have a Euclidean subgroup $[k,l]$ or $[l,m]$. The following theorem settles the case when $l=4$ or $6$.

\begin{THM}
\label{thmfourfour} 
Let $G \simeq [k,l,m]$ be crystallographic.  
Suppose the subgroup 
$[k,l]$ or $[l,m]$ is Euclidean, and that $l=4$ or $6$. Then $G^p$ is a C-group 
for any prime $p\geq 3$. 
\end{THM}

\noindent
\textbf {Proof.}  
Suppose $G_{3} = [k,l]$ (say) is Euclidean. The case $m=2$ was already settled, so let $m \neq 2$. For the moment, let $p$ be generic for $G$.
Then $V = \langle b_{0},\ldots,b_{3} \rangle$ is non-singular and $V_{3}$ is singular, 
so that every isometry of $V$ is uniquely determined by its  
effect  on $V_{3}$. Let $\mbox{\rm rad}\, V_{3} = \langle c \rangle$ (say).

Let $g \in G_{0}^p \cap G_{3}^p$.  Then $g$ leaves $V_{3}$ and $V_{0,3}$ invariant, and fixes $c$. Since $l=4$ or $6$, we necessarily have $G_{3} = [3,6]$ or $[4,4]$, so $G_{3}^p$ is 
a semi-direct product of its ``translation subgroup" $T^p$ (of order $p^2$) by $G_{0,3}^p$. 
In particular, since we are allowed to multiply $g$ by an element in $G_{0,3}^p$, we may 
assume that $g \in T^p$. We prove that this forces $g=e$, hence the desired conclusion.

When $G_{3} = [4,4]$, we may assume that the labels of the nodes $0$, $1$, $2$ of the diagram of $G$ are $1$, $2$, $1$, respectively. Then $c = b_{0} + b_{1} + b_{2}$ and
\[  T^{p} = \langle r_{0}r_{1}r_{2}r_{1}, r_{1}r_{0}r_{1}r_{2} \rangle 
\cong  \mathbb{Z}_{p} \times \mathbb{Z}_{p}. \]
Let $M^p$ denote the group of $p^2$ matrices of the form
\begin{equation}
\label{mlm}
M(\lambda,\mu) : =
\left[ \begin{array}{ccc}
1 & 0 & 0\\
0 & 1 & 0\\
\lambda&\mu&1
\end{array} \right] \qquad (\lambda,\mu \in \mathbb{Z}_p) .
\end{equation}
In the basis $b_{0},b_{1},c$ of $V_{3}$, each element of $T^{p}$, when restricted 
to $V_{3}$, is represented by a matrix in $M^p$. 
($T^p$ acts faithfully on $V_3$ by Lemma~\ref{3new1}.)
For example, we have
\[  (r_{0}r_{1}r_{2}r_{1})^{j} ( r_{1}r_{0}r_{1}r_{2})^k \mapsto M(2j,2(k-j)) \]
so $T^p$ and $M^p$ are clearly isomorphic.  Now suppose that $M(\lambda,\mu)$ is the matrix for $g$.  Then
\[ g(b_{1}) = b_{1} + \mu c = \mu b_{0} + (1+\mu) b_{1} + \mu b_{2}, \]
so we must have $\mu = 0$ because $V_{0,3}$ is invariant under $g$.  Similarly,
\[ \begin{array}{rll}
g(b_{2}) = g(c - b_{0} - b_{1})\!\!\! &= g(c) - g(b_{0}) - g(b_{1}) \\
\!\!\!&= c - (b_{0} + \lambda c) - b_{1} &\!\!\!\!\! = - \lambda b_{0} - 
\lambda b_{1} + (1 - \lambda)b_{2} , 
\end{array}\]
so also $\lambda = 0$, for the same reason. It follows that $g$ acts trivially on $V_{3}$ and 
hence also on $V$, that is, $g=e$.  

When $G_{3} = [3,6]$, we may take the labels of nodes $0$, $1$, $2$ of the diagram to be 
$1$, $1$, $3$, respectively. Then $c = b_{0} + 2b_{1} + b_{2}$ and 
\[  T^{p} = \langle r_{0}r_{1}(r_{2}r_{1})^{2}, r_{1}r_{0}(r_{1}r_{2})^{2} \rangle 
\cong  \mathbb{Z}_{p} \times \mathbb{Z}_{p}.  \]
The elements of $T^p$ still are represented by the matrices $M(\lambda,\mu)$ in $M^p$ (in the basis $b_{0},b_{1},c$ of $V_{3}$), so we can proceed in a similar fashion as above. In 
particular, if $M(\lambda,\mu)$ is the matrix for $g$, then we obtain $\mu = 0$ from 
\[ g(b_{1}) = b_{1} + \mu c = \mu b_{0} + (1+2\mu) b_{1} + \mu b_{2}, \]
and $\lambda = 0$ from 
\[  \begin{array}{rll}
g(b_{2}) = g(c - b_{0} - 2b_{1}) \!\!\! &= g(c) - g(b_{0}) - 2g(b_{1}) \\
\!\!\!& = c - (b_{0} + \lambda c) - 2b_{1} 
&\!\!\!\!\! = -\lambda b_{0} - 2\lambda b_{1} + (1-\lambda)b_{2} ,
\end{array}\]
in each case using the invariance of $V_{0,3}$ under $g$. Hence $g=e$, as required.

Only a few non-generic cases remain, which are not covered by
previous theorems: $p=3$ for $[4,4,6]$, $[3,6,3]$, $[3,6,4]$, $[3,6,6]$ or
$[3,6,\infty]$. But for each of the eleven possibly
distinct basic systems here, we easily verify the 
intersection condition  with the help of \textit{GAP}.
\hfill $\square$

This, then, leaves us with the groups $G=[6,3,m]$.  If $m=3$ or $4$, then $G_{0}$ is spherical, so Theorem~\ref{thmfourthree} applies and proves that $G^p$ is a C-group. It remains to investigate the cases $m=6$ or $\infty$. 

\begin{THM}
\label{thmfourfive} Let $G = [6,3,m]$ with $m=6$ or $\infty$.
Then for $m=6$,  $G^p$ is a C-group only for  $p=3$; and for $m =\infty$,
$G^p$ is a C-group  if and only if $p = 3$  
or $p \equiv \pm 5 \pmod{12}$.
\end{THM}

\noindent
\textbf{Proof.} As in  the previous theorem, we use \textit{GAP} to
verify the intersection condition in all cases with $p=3$. So suppose $p>3$.
We may assume that the nodes $0$, $1$, $2$, $3$ have the set of labels $3$, $1$, $1$, $3$, 
or $3$, $1$, $1$, $4$ according as $m=6$ or $\infty$. Since $V$ is non-singular, each isometry of $V$ is uniquely determined by its effect on the singular subspace $V_{3}$. In particular, 
\[ G_{3}^{p} \cong [6,3]^{p} \cong  T^p \rtimes \langle r_{0},r_{1} \rangle , \]
with
\[  T^{p} = \langle r_{1}r_{2}(r_{1}r_{0})^{2}, r_{2}r_{1}(r_{0}r_{1})^{2} \rangle 
\cong  \mathbb{Z}_{p} \times \mathbb{Z}_{p} .  \]
In the basis $b_{1},b_{2},c$ of $V_{3}$, with $c := b_{0} +2b_{1} + b_{2}$ (generating 
$\mbox{\rm rad}\, V_{3}$), each element of $T^{p}$, when restricted to $V_{3}$, is represented by a matrix $M(\lambda,\mu)$ as in (\ref{mlm}). In particular,
\[   r_{2}r_{1}(r_{0}r_{1})^{2} \mapsto M(-1,2),\quad  r_{1}r_{2}(r_{1}r_{0})^{2} 
	\mapsto M(1,1) . \]
Each element in $G_{3}^p$ is of the form $r_{1}^{i}(r_{0}r_{1})^{j} t$ with $i=0,1$, 
$j=0,\ldots,5$ and $t \in T^p$.  Inspection shows that the elements with $i=0$ have the following matrix representation in the basis $b_{1},b_{2},c$:\   if a matrix $M(a,b)$ represents an element $t=t(a,b)$ of $T^p$, then the matrix of $(r_{0}r_{1})^{j} t(a,b)$ is obtained from the matrix of $(r_{0}r_{1})^{j}$ by simply adding $a$ or $b$, respectively, to the first or second entry in its last row. 

Now suppose an element $g \in G_{3}^p$ leaves $V_{0,3}$ invariant.  Multiplying by $r_{1}$ if 
need be, we may assume that $i=0$, that is, $g = (r_{0}r_{1})^{j} t$ with $j = 0,\ldots,5$ and 
$t = t(a,b) \in T^p$.  The invariance of $V_{0,3}$ considerably restricts the possibilities for 
$j$ and $t$. In fact, inspection of the matrices for $(r_{0}r_{1})^{j} t(a,b)$ shows that we must 
have $(j,a,b) = (0,0,0)$, $(1,1,-1)$, $(2,1,-2)$, $(3,0,-2)$, $(4,-1,-1)$ or $(5,-1,0)$. 
Bearing in mind an initial multiplication by $r_{1}$, this leaves $12$ choices for $g$. 

In particular, since $G_{0}^p \cap G_{3}^p$ leaves $V_{0,3}$ invariant, it has order at most 
$12$ and contains $G_{0,3}^p$ as a subgroup of order $6$.  Comparison of the matrices shows that $(r_{0}r_{1})^{j} t(a,b)$ with $(j,a,b) = (1,1,-1)$, $(3,0,-2)$ or $(5,-1,0)$ is not 
contained in $G_{0,3}^p$.  It follows that $G^p$ fails to be a C-group if and only if one of these 
(and then all) three elements also lies in $G_{0}^p$.  Now consider 
\[ g := (r_{0}r_{1})^{3} t(0,-2) \] 
obtained for $(j,a,b)=(3,0,-2)$; its matrix in the basis $b_{1},b_{2},c$ is 
\begin{equation}
\label{mthree}
M =
\left[ \begin{array}{ccc}
-1 & 0 & 0\\
0 & -1 & 0\\
0 & 0  &1
\end{array} \right] .
\end{equation} 
Note that 
\begin{equation}
\label{gexp}
g =  (r_{0}r_{1})^{3}  (r_{1}r_{2}(r_{1}r_{0})^{2})^{\ell}  (r_{2}r_{1}(r_{0}r_{1})^{2})^{\ell} , 
\quad
3\ell \equiv -2 \bmod p .
\end{equation} 

First, if $G = [6,3,6]$, then $V_{0}$ is also singular, $d := b_{1} + 2b_{2} + b_{3}$ generates 
$\mbox{\rm rad}\, V_{0}$, and $b_{1},b_{2},d$ is a basis of $V_{0}$.  Now it is straightforward 
to check that $g(d)=d$; in fact, the equations 
\[ g(b_{3}) \cdot g(b_{j}) = b_{3} \cdot b_{j} \quad (j=0,\ldots,3)\] 
have the unique solution $g(b_{3}) = 2b_{1} + 4b_{2} + b_{3}$. We now exploit duality.  
In fact, by the symmetry of the diagram of $G$, we also have an element 
\[ g' = (r_{3}r_{2})^{3} t'(0,-2)\]  
(say) in $G_{0}^p$, which is obtained from $G_{0}^p$ in the same way as $g$ from 
$G_{3}^p$.  Since $g$ and $g'$ have the same $(4\times 4)$-matrices in the full basis 
$b_{1},b_{2},c,d$ of $V$, we must have $g=g'$.  Hence $g\in G_{0}^{p} \cap G_{3}^{p}$ but $g \not\in G_{0,3}^{p}$, so $G^p$ is not a C-group.  In particular, from (\ref{gexp}) and a similar expression for $g'$ we obtain the relation
\[ (r_{0}r_{1})^{3}  (r_{1}r_{2}(r_{1}r_{0})^{2})^{l}  (r_{2}r_{1}(r_{0}r_{1})^{2})^{l}  =
(r_{3}r_{2})^{3}  (r_{2}r_{1}(r_{2}r_{3})^{2})^{l}  (r_{1}r_{2}(r_{3}r_{2})^{2})^{l} , \]
with $3\ell \equiv -2 \bmod p$.

Finally, we consider  $G = [6,3,\infty]$. Here it is a little easier to 
work with $r := r_1 g$, still with $g$ defined in (\ref{gexp}). $G^p$ will be a C-group
exactly when $ r \not\in G_0^p$. But $r$ acts on $V_3$ as a reflection with root
$a:=b_1 + 2 b_2$; thus, since $V_3$ is a singular subspace of $V$, $r$ actually 
is a reflection in $G^p$. On the other hand, 
$G_0^p = O_1(V_0) \simeq O_1(3,p,0)$, for $p > 3$ (see \cite[5.7]{monsch}); and 
$a^2 = 3$. Hence, for $p>3$, $G^p$ is a C-group if and only if $3$ is a non-square modulo 
$p$, i.e. if and only if $p \equiv \pm 5 \pmod{12}$.
\hfill $\square$

We conclude this section with a look back at the
peculiar role of the non-generic prime $p = 3$, which is a frequent 
irritant in proofs but never an obstruction to polytopality.
We have seen for $p=3$ that when one of $k,l,m$ is $6$, different
basic systems for the crystallographic group $G = [k,l,m]$
can result in non-isomorphic reduced groups $G^3$. In all
cases, however, $G^3$ happens to be a C-group. Although this fact can be
checked by hand, we have   often resorted to 
computer verification.

In fact,  we can refine our description 
of the group $G^p$ in  singular cases.  
Indeed, if $n = \dim(V)$ and 
$r = \dim(\rm{rad}(V))$, then it  is easy  to see that
$$ \widehat{O}(V) \simeq \mathbb{Z}_p^{r(n-r)} 
\rtimes O(V/\mbox{rad}\;V)\;.
$$
In particular,   
if $n=4$ and $r=1,2 $ or $3$,  then it follows from the above isomorphism that (in non-generic cases),  $G^3$ must have order of the form $2^m 3^n$ (see \cite[p. 301]{monsch} for the orders of the groups of orthogonal type).  Several such groups appear in the
following sections. 
\vspace*{5mm}

\section{Groups $[k,l,m]$ with a spherical or Euclidean subgroup $[l,m]$}
\label{polfourone}

In the previous section we determined the crystallographic groups $G = [k,l,m]$ and primes 
$p \geq 3$ for which the modular reduction $G^p$ is C-group. We now investigate the corresponding regular $4$-polytope ${\cal P} = {\cal P}(G^p)$, whose automorphism group $\Gamma({\cal P})$ is $G^p$. When $p$ is generic for $G$, Lemma~\ref{3new1} implies that  
$G_{3}^p \simeq [k,l]^p$ and  $G_{0}^p \simeq [l,m]^p$. (In fact, this often holds
in non-generic situations, too.)  Thus the facets and vertex-figures  of ${\cal P}(G^p)$ 
are usually isomorphic to the regular maps associated with the reduced groups 
$[k,l]^p$ and $[l,m]^p$, respectively, as described
in \cite[\S\,5]{monsch}. These maps are orientable, 
because the even subgroups of $G_{0}^p$ and $G_{3}^p$ have index $2$.

The groups $G$ with disconnected diagrams $\Delta(G)$ were described early
in Section~\ref{interprop}, and the corresponding
polytopes are easily classified. Thus we shall assume from here on that 
$\Delta(G)$ is connected.  It then follows from Lemma 3.2 of 
\cite{monsch} that $G^p$ acts irreducibly on $V$, so long as
$\det(m_{ij}) \not\equiv 0 \pmod{p}$. (Recall that the $m_{ij}$ are the Cartan integers
appearing in \cite[\S 4]{monsch}. We take this opportunity to correct an
oversight in part (a) of  \cite[Lemma 3.2]{monsch}, where we should have stated  that the
roots $a_j$  form a basis for $V$. This actually is the case in all applications
here and in \cite{monsch}.) 

Before proceeding to specifics, we indicate how to decide whether $G^p$ is of orthogonal or
spherical type (see Section~\ref{basicreview}).  If, for example, $[k,l]$ is Euclidean, 
it is very easy  to check that a basic translation is a product of two 
reflections in  $[k,l]^p$  and has period $p$. By scanning the parameter $d(G)$ in 
Table 1 of \cite{monsch}, we see that the product of two reflections
has  period  at most $5$ in spherical cases. Thus
$G^p$ is of orthogonal type for $p>5$, and perhaps also for $p =3$ or $5$, cases which can be directly checked in \textit{GAP}. We will employ this
sort of analysis without much comment in what follows.

We break the discussion down into three cases according as the (vertex-figure) subgroup 
$[l,m]$ of $G$ is spherical, Euclidean or hyperbolic, respectively. In this section we treat the groups $G$ with a spherical or Euclidean subgroup $[l,m]$; groups with hyperbolic subgroups $[l,m]$ are studied in the next section. We begin with the spherical case.
\medskip

\noindent
{\bf 4.1.\ Groups with a spherical subgroup $[l,m]$}

Let $G = [k,l,m]$ be crystallographic, let $[l,m]$ be spherical, and let $p\geq 3$. Then 
$G^p$ is a C-group by Theorem~\ref{thmfourthree}.  
Moreover, $[l,m]^p \cong [l,m]$, so the vertex-figures of the corresponding 
polytope $\cal P$ are isomorphic to Platonic solids $\{l,m\}$.  The finite or Euclidean groups 
$G$ were already discussed in \cite[\S 5-6]{monsch}, so we may assume here that $[k,l]$ is not spherical. 

If $[k,l]$ is Euclidean, then $G=[4,4,3]$, $[6,3,3]$ or $[6,3,4]$, and $\cal P$ is locally 
toroidal. Recall that a regular $4$-polytope is {\em locally toroidal\/} if its facets and 
vertex-figures are toroidal or spherical, but not all spherical. 

For $G =[4,4,3]$ with diagram
\[ \stackrel{1}{\bullet}\!\frac{}{\;\;\;\;\;\;}\!
\stackrel{2}{\bullet}\!\frac{}{\;\;\;\;\;\;}\!
\stackrel{1}{\bullet}\!\frac{}{\;\;\;\;\;\;}\!
\stackrel{1}{\bullet} ,\]
the facets of $\cal P$ are toroidal maps $\{4,4\}_{(p,0)}$ and the vertex-figures are $3$-cubes 
$\{4,3\}$. Now $\mbox{disc}(V) \sim -1$, independent of $p$, so $\epsilon(V) = 1$ if and only if 
$p \equiv 1 \bmod 4$. Moreover, $G^p = O_{1}(V)$ or $O(V)$ according as $2$ is a square or 
non-square, that is, $p \equiv \pm 1 \bmod 8$ or  $p \equiv \pm 3 \bmod 8$. Hence $\cal P$ has automorphism group
\begin{equation}\label{443ways}
 \Gamma({\cal P}) = G^p = 
\left\{ \begin{array}{ll}
O_1(4,p,1)\;, & \mbox{ if } p \equiv 1 \pmod{8}\\
O_1(4,p,-1)\;, & \mbox{ if } p \equiv 7 \pmod{8}\\
O(4,p,1)\;, & \mbox{ if } p \equiv  5 \pmod{8}\\
O(4,p,-1)\;, & \mbox{ if } p \equiv  3 \pmod{8}.
\end{array}\right. 
\end{equation}
In particular, for $p=3$ we obtain the group 
\begin{equation}
\label{unfofoth}
\Gamma({\cal P}) \,=\, O(4,3,-1)\, \cong \, O_{1}(4,3,-1) \times C_{2}\, \cong\, 
S_{6}\times C_{2} 
\end{equation} 
(\cite[\S 3]{monsch}). In this case,
\begin{equation}
\label{fofothreezero}
{\cal P} = \{\{4,4\}_{(3,0)},\{4,3\}\} ,
\end{equation}
the universal regular $4$-polytope with facets $\{4,4\}_{(3,0)}$ and vertex-figures $\{4,3\}$ 
(see \cite[Thm. 10B3]{arp}, which also implies that the product in (\ref{unfofoth}) indeed is 
direct).  Recall that $\{{\cal P}_{1},{\cal P}_{2}\}$ denotes the universal regular 
$(n+1)$-polytope (if it exists) with facets isomorphic to ${\cal P}_{1}$ and vertex-figures 
isomorphic to ${\cal P}_{2}$ (see \cite[Ch. 4]{arp}).

For $G = [6,3,m]$, with $m=3$ or $4$, we take the diagram
\begin{equation}
\label{diastm}
\stackrel{3}{\bullet}\!\frac{}{\;\;\;\;\;\;}\!
\stackrel{1}{\bullet}\!\frac{}{\;\;\;\;\;\;}\!
\stackrel{1}{\bullet}\!\frac{}{\;\;\;\;\;\;}\!
\stackrel{1}{\bullet} 
\qquad \mbox{or} \qquad
\stackrel{3}{\bullet}\!\frac{}{\;\;\;\;\;\;}\!
\stackrel{1}{\bullet}\!\frac{}{\;\;\;\;\;\;}\!
\stackrel{1}{\bullet}\!\frac{}{\;\;\;\;\;\;}\!
\stackrel{2}{\bullet} ,
\end{equation}
respectively. Now the vertex-figures of $\cal P$ are tetrahedra $\{3,3\}$ or octahedra 
$\{3,4\}$, respectively, and the facets are toroidal maps $\{6,3\}_{(p,0)}$ 
for all $p\geq 3$ (see \cite[\S 5.6]{monsch}).  Note that 
$\mbox{disc}(V) \sim -3$, so $V$ is non-singular if $p>3$, and $\epsilon(V) = 1$ if and only if 
$-3$ is a square; furthermore,
$G^p = O_{1}(V)$ or $O(V)$ depending on the quadratic character of $3$
(or $2$ and $3$).

In particular, if $G = [6,3,3]$ and $p>3$, we have 
\begin{equation}
\label{grsithth}
\Gamma({\cal P}) = 
\left\{ \begin{array}{ll}
O_1(4,p,1)\;, & \mbox{ if } p \equiv 1 \pmod{12}\\
O_1(4,p,-1)\;, & \mbox{ if } p \equiv 11 \pmod{12}\\
O(4,p,1)\;, & \mbox{ if } p \equiv 7 \pmod{12}\\
O(4,p,-1)\;, & \mbox{ if } p \equiv 5 \pmod{12}.
\end{array}\right. 
\end{equation}
When $p=3$ we also obtain a C-group $G^3$, which acts on a singular 
space $V$ and has order $1296$. The corresponding subgroup $G_{3}^3$ has order $108$,
and the facets are toroidal maps $\{6,3\}_{(3,0)}$ (see \cite[\S\,5.6]{monsch}).
Note that Lemma~\ref{3new1} does not apply; indeed the subdiagram on nodes
$0,1,2$ of $\Delta(G)$ defines the rank $3$ group $H^3$, of order just $36$, for the
map $\{6,3\}_{(1,1)}$. In a sense, the subspace $V_3$ of $V$ cannot fully represent
the structure of the facet. The alternative diagram
\[ \stackrel{1}{\bullet}\!\frac{}{\;\;\;\;\;\;}\!
\stackrel{3}{\bullet}\!\frac{}{\;\;\;\;\;\;}\!
\stackrel{3}{\bullet}\!\frac{}{\;\;\;\;\;\;}\!
\stackrel{3}{\bullet}\]
represents the same group, but without this deficiency. 
By comparing group orders, we find that we have the universal polytope
\begin{equation}
\label{sithun}
{\cal P} = \{\{6,3\}_{(3,0)},\{3,3\}\} 
\end{equation}
with
\[ G^{3} = \Gamma({\cal P}) \cong {[1\,1\,2]}^{3}\rtimes C_{2} , \]
where ${[1\,1\,2]}^{3}$ is a certain unitary reflection group in $\mathbb{C}^4$ 
(see \cite[Thm. 11B5]{arp}).  (The superscript on ${[1\,1\,2]}^{3}$ signifies a group relation, 
not reduction modulo $3$).

When $G = [6,3,4]$ and $p>3$, we similarly find that
\begin{equation}\label{case634} 
\Gamma({\cal P}) = 
\left\{ \begin{array}{ll}
O_1(4,p,1)\;, & \mbox{ if } p \equiv 1 \pmod{24}\\
O_1(4,p,-1)\;, & \mbox{ if } p \equiv 23 \pmod{24}\\
O(4,p,1)\;, & \mbox{ if } p \equiv 7,13,19 \pmod{24}\\
O(4,p,-1)\;, & \mbox{ if } p \equiv 5,11,17 \pmod{24}.
\end{array}\right. \end{equation}
For $p=3$ and for either of the diagrams
\begin{equation}\label{634mod3}
 \stackrel{3}{\bullet}\!\frac{}{\;\;\;\;\;\;}\!
\stackrel{1}{\bullet}\!\frac{}{\;\;\;\;\;\;}\!
\stackrel{1}{\bullet}\!\frac{}{\;\;\;\;\;\;}\!
\stackrel{2}{\bullet}
\qquad \mbox{ or }\qquad
\stackrel{1}{\bullet}\!\frac{}{\;\;\;\;\;\;}\!
\stackrel{3}{\bullet}\!\frac{}{\;\;\;\;\;\;}\!
\stackrel{3}{\bullet}\!\frac{}{\;\;\;\;\;\;}\!
\stackrel{6}{\bullet}\;\;, \end{equation}
 $G^3$ has order $2592$;  we obtain the same polytope
$\mathcal{P}$ of type $\{\{6,3\}_{(3,0)},\{3,4\}\}$.
However, $\mathcal{P}$ is not the (infinite!)  universal polytope of this type 
\cite[Thm. 11B5]{arp}.

We now consider the case that the subgroup $[k,l]$ is hyperbolic 
(and $[l,m]$ is still spherical). 
The corresponding groups are $G=[\infty,3,3]$, $[\infty,3,4]$, $[6,4,3]$ or $[\infty,4,3]$. 

For $G=[\infty,3,m]$, with $m=3$ or $4$, we employ the diagram 
\[ \stackrel{4}{\bullet}\!\frac{}{\;\;\;\;\;\;}\!
\stackrel{1}{\bullet}\!\frac{}{\;\;\;\;\;\;}\!
\stackrel{1}{\bullet}\!\frac{}{\;\;\;\;\;\;}\!
\stackrel{1}{\bullet}
\qquad \mbox{ or }\qquad
\stackrel{4}{\bullet}\!\frac{}{\;\;\;\;\;\;}\!
\stackrel{1}{\bullet}\!\frac{}{\;\;\;\;\;\;}\!
\stackrel{1}{\bullet}\!\frac{}{\;\;\;\;\;\;}\!
\stackrel{2}{\bullet}, \]
respectively. In either case,
 $G_{3}^{p}=[\infty,3]^{p}$, so that the facets of $\cal P$ are the 
regular maps $\mathcal{M}_{p,3}$
of type $\{p,3\}$ described in \cite[\S 5.7]{monsch} (see also \cite{regmaps}); 
the vertex-figures 
are tetrahedra $\{3,3\}$ or octahedra $\{3,4\}$, respectively. In particular, if $p=3$, 
we obtain the $4$-simplex $\{3,3,3\}$ or 4-cross-polytope $\{3,3,4\}$, respectively.  Now 
$\mbox{disc}(V) \sim -1$ or $-2$, respectively, independent of $p$; and $G^p = O_{1}(V)$ unless $m = p = 3$, or $m=4$ and $2$ is a non-square. Thus, when $m=3$ we have
\begin{equation} \label{caseinf33}
\Gamma({\cal P}) = 
\left\{ \begin{array}{ll}
S_5\;, & \mbox{ if }  p =3\\
O_1(4,p,1)\;, & \mbox{ if }  p \equiv 1 \pmod{4}\\
O_1(4,p,-1)\;, & \mbox{ if } p \equiv 3 \pmod{4} , \; (p \neq 3),
\end{array}\right. \end{equation}
and when $m=4$ and $p>3$, we have
\begin{equation}\label{caseinf34}
\Gamma({\cal P}) = 
\left\{ \begin{array}{ll}
O_1(4,p,1)\;, & \mbox{ if } p \equiv 1 \pmod{8}\\
O_1(4,p,-1)\;, & \mbox{ if } p \equiv 7 \pmod{8}\\
O(4,p,1)\;, & \mbox{ if } p \equiv 3 \pmod{8}\\
O(4,p,-1)\;, & \mbox{ if } p \equiv 5 \pmod{8} .
\end{array}\right. \end{equation}
(When $m=4$ and $p=3$, the group is $B_4$, which is of index $3$ in
$O(4,3,1) \cong F_4$.) For both $m=3$ and $m=4$, and all $p \geq 3$, 
we have  $G_{3}^{p} = O_{1}(3,p,0)$, of order $p(p^{2}-1)$. 
For $p=5$   the facets are 
isomorphic to   dodecahedra $\{5,3\}  =  \mathcal{M}_{5,3}$; 
and for $p=7$ they are isomorphic to Klein's map 
$\{7,3\}_8  = \mathcal{M}_{7,3}$, of genus $3$ (see \cite[\S\,8.6]{genrel}). 

With $m=3$, we here encounter, for the first time, the  
classical regular $4$-polytopes of `pentagonal type'.  For future reference, let us 
denote $\mathcal{P}$ by $\mathcal{C}_{p,3,3}$\,. Likewise, as suggested above,
it suits us to denote the regular maps arising from

\[ \stackrel{4}{\bullet}\!\frac{}{\;\;\;\;\;\;}\!
\stackrel{1}{\bullet}\!\frac{}{\;\;\;\;\;\;}\!
\stackrel{1}{\bullet}
\qquad \mbox{ and }\qquad
\stackrel{4}{\bullet}\!\frac{}{\;\;\;\;\;\;}\!
\stackrel{1}{\bullet}\!\frac{}{\;\;\;\;\;\;}\!
\stackrel{4}{\bullet} \]
by $\mathcal{M}_{p,3}$ and $\mathcal{M}_{p,p}$\,, respectively (see \cite[\S 5.7]{monsch}).

Thus,  for  $p=5$, $\mathcal{C}_{5,3,3}$  is  
the $120$-cell $\{5,3,3\}$ (isomorphic to both the convex regular polytope
of this type and to the star-polytope $\{\frac{5}{2},3,3\}$; see \cite[7D]{arp}). 
For $p=7$, we find that $\mathcal{C}_{7,3,3}$  is 
the universal regular polytope
\[ \{\{7,3\}_{8},\{3,3\}\} , \]
first described in \cite{monw1}.  Rephrasing the conjectured presentation
in \cite[3.1]{monw1}, we now offer  

\noindent
{\bf Conjecture~1}: For primes $p \geq 3$, $\mathcal{C}_{p,3,3}$ is the universal
polytope  
$$
\{ \, \mathcal{M}_{p,3}\, , \, \{3,3\}\, \}\;\;.
$$
(This has been verified by coset enumeration for $p \leq 31$. We cannot, in general, 
expect the facets to be determined in some elegant way, say by  their Petrie polygons or holes.)
\medskip

Moving on to the next class of groups,
let $G=[k,4,3]$, with $k=6$ or $\infty$, and let the diagram be
\[ \stackrel{3}{\bullet}\!\frac{}{\;\;\;\;\;\;}\!
\stackrel{1}{\bullet}\!\frac{}{\;\;\;\;\;\;}\!
\stackrel{2}{\bullet}\!\frac{}{\;\;\;\;\;\;}\!
\stackrel{2}{\bullet}
\qquad \mbox{ or }\qquad
\stackrel{4}{\bullet}\!\frac{}{\;\;\;\;\;\;}\!
\stackrel{1}{\bullet}\!\frac{}{\;\;\;\;\;\;}\!
\stackrel{2}{\bullet}\!\frac{}{\;\;\;\;\;\;}\!
\stackrel{2}{\bullet}, \]
respectively.  Then $G_{3}^{p}=[k,4]^{p}$, so the facets are   maps of type $\{6,4\}$ or 
$\{p,4\}$, respectively; see \cite[\S 5.8-5.9]{monsch} for details. 
Of course, the vertex-figures are cubes $\{4,3\}$. 
Since $\mbox{disc}(V) \sim -15$ or 
$-2$, respectively,   $V$ is singular if $k=6$ and $p=3$ or $5$.

Suppose $k=6$.
For $p=3$, the group $G^3$ has order $2592$ and appeared earlier in (\ref{634mod3}) 
with different generators. In this new guise, $G^3$ is the group of the
universal polytope of type  
\begin{equation}\label{univ64sub43}
\{ \{6,4\}_4 \, ,  \, \{4,3\} \} \;\;.
\end{equation}
Indeed, for $p=3$, each basic system  gives just this   polytope,
whose facets are isomorphic to the map $\{6,4\}_4$ (the Petrial of the toroidal map 
$\{4,4\}_{(3,3)}$). 
(See \cite[7B2]{arp} for a description
of $\mathcal{Q}^\pi$, the \textit{Petrial} or
\textit{Petrie dual} of a map $\mathcal{Q}$.)

For $p=5$, the facets
are Coxeter-Petrie polyhedra $\{6,4\!\mid\! 3\}$ 
(see \cite{genrel}).  Since $G^5$ has order $30000$, the corresponding polytope has 
$125$ of these facets and $625$ vertices.
  
Otherwise, if $G = [6,4,3]$ and $p > 5$, then we have
\[ \Gamma({\cal P}) = 
\left\{ \begin{array}{ll}
O_1(4,p,1)\;,  & \mbox{ if }  p \equiv 1,23,47,49  \pmod{120}\\
O_1(4,p,-1)\;, & \mbox{ if }  p \equiv 71,73,97,119 \pmod{120}\\
O(4,p,1)\;,  & \mbox{ if }  p \equiv 17,19,31,53,61,77,79,83,91,107,109,113 \pmod{120}\\
O(4,p,-1)\;, & \mbox{ if }  p \equiv 7,11,13,29,37,41,43,59,67,89,101,103 \pmod{120}.
\end{array}\right. \]

Finally, when $G = [\infty,4,3]$, the polytopes are of type  
$\{p,4,3\}$ and their groups are 
those described  in (\ref{caseinf34}) above (now allowing $p=3$), with new generators of course. For $p=3$ we obtain the $24$-cell $\{3,4,3\}$. When $p=5$, the polytope has facets isomorphic to Gordan's map $\{5,4\}_{6}$ of genus $4$ (see \cite{genrel}). 
\medskip

\noindent
{\bf 4.2.\ Groups with a Euclidean subgroup $[l,m]$}

Next we consider the crystallographic groups $G = [k,l,m]$ with a Euclidean subgroup $[l,m]$. The groups with a spherical subgroup $[k,l]$ have already been discussed in the dual setting, so we may further assume that $[k,l]$ is Euclidean or hyperbolic. By Theorems~\ref{thmfourfour} and \ref{thmfourfive}, $G^p$ is a C-group for all $p\geq 3$, with these exceptions:\  for  $G=[6,3,6]$ only $p=3$ is acceptible, and for $G=[\infty,3,6]$, only
$p =3$ and $p \equiv \pm5 \pmod{12}$. 
The corresponding regular $4$-polytopes 
$\cal P$ all have toroidal vertex-figures. 

For the group $G = [k,4,4]$, with $k=4$, $6$ or $\infty$, we take the diagram 
\[ \stackrel{2}{\bullet}\!\frac{}{\;\;\;\;\;\;}\!
\stackrel{1}{\bullet}\!\frac{}{\;\;\;\;\;\;}\!
\stackrel{2}{\bullet}\!\frac{}{\;\;\;\;\;\;}\!
\stackrel{1}{\bullet} \,,\;
\qquad 
\stackrel{3}{\bullet}\!\frac{}{\;\;\;\;\;\;}\!
\stackrel{1}{\bullet}\!\frac{}{\;\;\;\;\;\;}\!
\stackrel{2}{\bullet}\!\frac{}{\;\;\;\;\;\;}\!
\stackrel{1}{\bullet}
\qquad \mbox{ or }\qquad
\stackrel{4}{\bullet}\!\frac{}{\;\;\;\;\;\;}\!
\stackrel{1}{\bullet}\!\frac{}{\;\;\;\;\;\;}\!
\stackrel{2}{\bullet}\!\frac{}{\;\;\;\;\;\;}\!
\stackrel{1}{\bullet}, \]
respectively. Their polytopes $\cal P$ have toroidal vertex-figures $\{4,4\}_{(p,0)}$ and 
facets isomorphic to the maps of type $\{k,4\}$ described in \cite[\S 5.5,5.8,5.9]{monsch}. Now 
$\mbox{disc}(V) \sim -1$ in each case, except when $k=6$ and $p=3$; in this latter case, $V$ is singular (for each admissible diagram).  

When $G = [4,4,4]$, the polytope is self-dual and its facets are also maps $\{4,4\}_{(p,0)}$. The group $\Gamma({\cal P})$ is the same as in (\ref{443ways}); that is, we have 
\[ [4,4,4]^{p} \cong [4,4,3]^{p} \]
for each $p\geq 3$. In fact, the Coxeter group $[4,4,4]$ is known to be a subgroup of index 
$3$ in $[4,4,3]$, and under the modular reduction this index collapses to $1$; see, for example, \cite[\S 10E]{arp}, which also explains the corresponding relationship between the polytopes. In particular, if $[4,4,3] = \langle r_{0},\ldots,r_{3} \rangle$ (say), then we can identify $[4,4,4]$ with the subgroup 
$\langle r_{1},r_{0},r_{2}r_{1}r_{2},r_{3} \rangle$, with the generators taken in this order; 
modulo $p$, this subgroup is the group itself.  

For $p=3$ we obtain the universal regular polytope
\[ \{\{4,4\}_{(3,0)},\{4,4\}_{(3,0)}\} \]
with group $S_{6} \times C_{2}$, which is related to the polytope in (\ref{fofothreezero}) 
(see \cite[10E6 and Thm. 10C12]{arp}).

For $G = [6,4,4]$ and $p>3$, we have
\[ \Gamma({\cal P}) = 
\left\{ \begin{array}{ll}
O_1(4,p,1)\;,  & \mbox{ if }  p \equiv 1 \pmod{24}\\
O_1(4,p,-1)\;, & \mbox{ if }  p \equiv 23 \pmod{24}\\
O(4,p,1)\;,  & \mbox{ if }  p \equiv 5,13,17 \pmod{24}\\
O(4,p,-1)\;, & \mbox{ if }  p \equiv 7,11,19 \pmod{24}.
\end{array}\right. \]
When $p=5$ the facets of $\mathcal{P}$ are Coxeter-Petrie polyhedra 
$\{6,4\!\mid\! 3\}$ (see \cite{genrel}). For 
$p=3$,  the facets are maps $\{6,4\}_4$,   the vertex-figures are maps of type
$\{4,4\}_{(3,0)}$, and the  diagram given above yields a group $G^3$ with order $1296$.
However, the alternate diagram
$$\stackrel{1}{\bullet}\!\frac{}{\;\;\;\;\;\;}\!
\stackrel{3}{\bullet}\!\frac{}{\;\;\;\;\;\;}\!
\stackrel{6}{\bullet}\!\frac{}{\;\;\;\;\;\;}\!
\stackrel{3}{\bullet} $$ 
results instead in a group $G^3$ of order $3888$. We thus obtain different polytopes with
the same local structure.

When $G = [\infty,4,4]$, the polytopes are of type $\{p,4,4\}$, and although
differently generated,  the  group $\Gamma(\mathcal{P})$ is again described by
(\ref{443ways}) when $p\geq 3$.
For $p=3$ we obtain the universal regular polytope 
\[ \{\{3,4\},\{4,4\}_{(3,0)}\} , \]
the dual of (\ref{fofothreezero}), with group $S_{6} \times C_{2}$. For $p=5$ the polytope 
$\cal P$ has facets isomorphic to Gordan's map $\{5,4\}_{6}$ of genus $4$. 
\medskip

Next we investigate the groups $[k,3,6]$.  By Theorem~\ref{thmfourfive},
when   $k=6$   we need only consider $p=3$. The diagrams 

\begin{equation}\label{636othera}
\stackrel{1}{\bullet}\!\frac{}{\;\;\;\;\;\;}\!
\stackrel{3}{\bullet}\!\frac{}{\;\;\;\;\;\;}\!
\stackrel{3}{\bullet}\!\frac{}{\;\;\;\;\;\;}\!
\stackrel{1}{\bullet} \;\;\mbox{ and } \;\;
\stackrel{1}{\bullet}\!\frac{}{\;\;\;\;\;\;}\!
\stackrel{3}{\bullet}\!\frac{}{\;\;\;\;\;\;}\!
\stackrel{3}{\bullet}\!\frac{}{\;\;\;\;\;\;}\!
\stackrel{9}{\bullet} \end{equation}
give isomorphic groups $G^3$ of order $1944$, yet \textit{distinct} polytopes of type
$\{ \{6,3\}_{(3,0)} , \{3,6\}_{(3,0)}  \}$ (with the same type of facets and 
vertex-figures, the latter type  being the dual of the former). 
However, only the first of these two $4$-polytopes is in itself self-dual. The other
essentially distinct diagram

\begin{equation}\label{636otherb}
\stackrel{3}{\bullet}\!\frac{}{\;\;\;\;\;\;}\!
\stackrel{1}{\bullet}\!\frac{}{\;\;\;\;\;\;}\!
\stackrel{1}{\bullet}\!\frac{}{\;\;\;\;\;\;}\!
\stackrel{3}{\bullet} \end{equation}
yields a group $G^3$ of order $216$ and the universal self-dual
polytope of type

\begin{equation}\label{636univ}
\{ \{6,3\}_{(1,1)} , \{3,6\}_{(1,1)}  \} 
\end{equation}
(see \cite[11C8]{arp}).
\medskip

Let us now turn to
the case $k=\infty$, with diagram
\begin{equation}
\label{diakts}
\stackrel{4}{\bullet}\!\frac{}{\;\;\;\;\;\;}\!
\stackrel{1}{\bullet}\!\frac{}{\;\;\;\;\;\;}\!
\stackrel{1}{\bullet}\!\frac{}{\;\;\;\;\;\;}\!
\stackrel{3}{\bullet} . 
\end{equation}
Now $\mbox{disc}(V) \sim -3$, so $V$ is singular if $p=3$. When $p>3$, the vertex-figures of 
$\cal P$ are toroidal maps $\{3,6\}_{(p,0)}$ (see \cite[\S 5.6]{monsch}), and the facets are 
the maps of type $\{p,3\}$ described in \cite[\S 5.7]{monsch}. 
Keeping in mind Theorem~\ref{thmfourfive}, we have just two possibilities:

\begin{equation}
\label{infthsix}
\Gamma({\cal P}) = 
\left\{ \begin{array}{ll}
O(4,p,1)\;,  & \mbox{ if }  p \equiv 7 \pmod{12}\\
O(4,p,-1)\;, & \mbox{ if }  p \equiv 5 \pmod{12}.
\end{array}\right. 
\end{equation}
For $p=5$ the facets of $\cal P$ are dodecahedra $\{5,3\}$, and for $p=7$ they are isomorphic to 
Klein's map $\{7,3\}_8$ of genus $3$.  When $p=3$, both the diagram in
(\ref{diakts}) and its alternative
yield the universal regular polytope 
\[  {\cal P} = \{\{3,3\},\{3,6\}_{(3,0)}\}, \]
namely the dual of (\ref{sithun}). Again the group  is ${[1\,1\,2]}^{3} \rtimes C_{2}$ of order $1296$. 
\medskip

It remains to study the groups $G = [k,6,3]$ with $k=3$, $4$, $6$ or $\infty$, where we take the diagrams 
\[ \stackrel{1}{\bullet}\!\frac{}{\;\;\;\;\;}\!
\stackrel{1}{\bullet}\!\frac{}{\;\;\;\;\;}\!
\stackrel{3}{\bullet}\!\frac{}{\;\;\;\;\;}\!
\stackrel{3}{\bullet}\,,
\quad 
\stackrel{2}{\bullet}\!\frac{}{\;\;\;\;\;}\!
\stackrel{1}{\bullet}\!\frac{}{\;\;\;\;\;}\!
\stackrel{3}{\bullet}\!\frac{}{\;\;\;\;\;}\!
\stackrel{3}{\bullet}\,,
\quad 
\stackrel{3}{\bullet}\!\frac{}{\;\;\;\;\;}\!
\stackrel{1}{\bullet}\!\frac{}{\;\;\;\;\;}\!
\stackrel{3}{\bullet}\!\frac{}{\;\;\;\;\;}\!
\stackrel{3}{\bullet}
\quad \mbox{ or }\quad
\stackrel{4}{\bullet}\!\frac{}{\;\;\;\;\;}\!
\stackrel{1}{\bullet}\!\frac{}{\;\;\;\;\;}\!
\stackrel{3}{\bullet}\!\frac{}{\;\;\;\;\;}\!
\stackrel{3}{\bullet}, \]
respectively. When $p>3$ the vertex-figures of $\cal P$ are 
toroidal maps $\{6,3\}_{(p,0)}$, and the 
facets are the maps of type $\{k,6\}$ or $\{p,6\}$, for  $k=3,4,6$ or $k=\infty$, as described in 
\cite[\S 5.6,5.8,5.10,5.11]{monsch}.  Since  $\mbox{disc}(V) \sim -3$ in each case,   $V$ is 
 singular if $p=3$. 

When $G = [3,6,3]$ and $p>3$, the polytope $\cal P$ is self-dual and its facets are toroidal maps $\{3,6\}_{(p,0)}$. The group $\Gamma({\cal P})$ is 
the same as in (\ref{grsithth}); that is, we have 
\begin{equation}
\label{subrel}
[3,6,3]^{p} \cong [6,3,3]^{p} 
\end{equation}
for $p>3$. In fact, the Coxeter group $[3,6,3]$ is a subgroup of index $4$ in $[6,3,3]$, and 
under reduction modulo $p$ this index becomes $1$, again so long as $p>3$; 
see, for example, \cite[Sect. 11G,11H]{arp}, 
which also describes the relationship between the polytopes. In particular, if 
$[6,3,3] = \langle r_{0},\ldots,r_{3} \rangle$ (say), then $[3,6,3]$ can be identified with the 
subgroup $\langle r_{0},r_{1}r_{0}r_{1},r_{2},r_{3} \rangle$, with the generators taken in this 
order; modulo $p$, this is the whole group. 

When $p=3$, we obtain the (combinatorially flat) universal regular polytope
\begin{equation}\label{363sub3univ} 
{\cal P} = \{\{3,6\}_{(1,1)},\{6,3\}_{(3,0)}\}
\end{equation}   
with group
\[ \Gamma({\cal P}) \cong {[1\,1\,1]}^{3} \rtimes S_{3} ,\]
of order $324$, where ${[1\,1\,1]}^{3}$ denotes a certain unitary reflection group in 
$\mathbb{C}^3$ (see \cite[Thm. 11E7]{arp}). 
Of course, $\mathcal{P}$ is not self-dual, although the alternative diagram 
\[ \stackrel{3}{\bullet}\!\frac{}{\;\;\;\;\;\;}\!
\stackrel{3}{\bullet}\!\frac{}{\;\;\;\;\;\;}\!
\stackrel{1}{\bullet}\!\frac{}{\;\;\;\;\;\;}\!
\stackrel{1}{\bullet} \]
does yield  the dual. Note also that the isomorphism in (\ref{subrel}) must  
fail for $p=3$, regardless of choice of diagrams. Indeed, $[3,6,3]^3$   does have index $4$
in $[6,3,3]^3$.

Incidentally, we also have 
\[ [6,3,6]^{p} \cong [6,3,3]^{p}  \]
for $p>3$; but then, as we have seen, 
$[6,3,6]^p$ is not a C-group with its natural generators. (In this context, note that 
\cite[Thm. 11H10]{arp} is incorrect for the parameter vectors $(s,0)$ with $s$ not divisible by 
$3$, so  does not yield any polytopes of type $\{6,3,6\}$.)

If $G = [4,6,3]$ and $p>3$, then $\Gamma(\mathcal{P})$ is the group
described earlier in (\ref{case634}), though again with new generators.  
For $p=5$ the facets of $\cal P$ are Coxeter-Petrie polyhedra 
$\{4,6\!\mid\! 3\}$ (see \cite{genrel}).  For $p=3$, the facets are   
maps $\{4,6\}_4$; but the diagrams
$$
\stackrel{2}{\bullet}\!\frac{}{\;\;\;\;\;}\!
\stackrel{1}{\bullet}\!\frac{}{\;\;\;\;\;}\!
\stackrel{3}{\bullet}\!\frac{}{\;\;\;\;\;}\!
\stackrel{3}{\bullet}\;\; , \;\;
\stackrel{6}{\bullet}\!\frac{}{\;\;\;\;\;}\!
\stackrel{3}{\bullet}\!\frac{}{\;\;\;\;\;}\!
\stackrel{1}{\bullet}\!\frac{}{\;\;\;\;\;}\!
\stackrel{1}{\bullet}$$
lead respectively to vertex-figures $\{6,3\}_{(3,0)}$ and group order $3888$,
or to vertex-figures $\{6,3\}_{(1,1)}$ and group order $432$.

When $G = [6,6,3]$ and $p>3$, the facets of $\cal P$ are self-dual maps of type $\{6,6\}$, with 
group $O_{1}(3,p,0)$ or $O(3,p,0)$ according as $p \equiv \pm 1 \bmod{12}$ or 
$p \not\equiv \pm 1 \bmod{12}$. In fact, $\Gamma({\cal P})$ is   
a newly generated version of  the group described in (\ref{grsithth}).
When $p=3$, we find that the type of facets depends on the diagram chosen; they are isomorphic to the Petrials of two maps first described by Sherk in \cite{sherk}\,: \
$\{6,6\}_{(1,1)}$, with group of order $72$ and genus $4$, and
$\{6,6\}_{(3,0)}$, with group $216$ and genus $10$. Note that 
$\{6,6\}_{(1,1)}  \cong \{6,6 \mid 2 \}$ (see \cite[\S\, 8.5]{genrel}).
The various diagrams yield 
four distinct universal polytopes, summarized in the chart below: 

\begin{equation}\label{univ663}
\begin{array}{c|c|c} 
\Delta(G)  & | G^3 |  &  \mbox{ The  Universal Polytope  }\\
 & &  \\\hline
  & &  \\
\stackrel{1}{\bullet}\!\frac{}{\;\;\;\;\;}\!
\stackrel{3}{\bullet}\!\frac{}{\;\;\;\;\;}\!
\stackrel{1}{\bullet}\!\frac{}{\;\;\;\;\;}\!
\stackrel{1}{\bullet}&  216& \{ \{6,6\}^{\pi}_{(1,1)} \, , \,  \{6,3\}_{(1,1)} \}\\
 & &  \\ \hline
  & &  \\
\stackrel{9}{\bullet}\!\frac{}{\;\;\;\;\;}\!
\stackrel{3}{\bullet}\!\frac{}{\;\;\;\;\;}\!
\stackrel{1}{\bullet}\!\frac{}{\;\;\;\;\;}\!
\stackrel{1}{\bullet}& 648 & \{ \{6,6\}^{\pi}_{(3,0)} \, , \,   \{6,3\}_{(1,1)} \} \\
 & &  \\ \hline 
  & &  \\
\stackrel{3}{\bullet}\!\frac{}{\;\;\;\;\;}\!
\stackrel{1}{\bullet}\!\frac{}{\;\;\;\;\;}\!
\stackrel{3}{\bullet}\!\frac{}{\;\;\;\;\;}\!
\stackrel{3}{\bullet}& 648 & \{ \{6,6\}^{\pi}_{(1,1)} \, , \, \{6,3\}_{(3,0)} \}\\
 & &  \\ \hline
  & &  \\
\stackrel{1}{\bullet}\!\frac{}{\;\;\;\;\;}\!
\stackrel{3}{\bullet}\!\frac{}{\;\;\;\;\;}\!
\stackrel{9}{\bullet}\!\frac{}{\;\;\;\;\;}\!
\stackrel{9}{\bullet} & 5832  & \{ \{6,6\}^{\pi}_{(3,0)}\, , \, \{6,3\}_{(3,0)} \}\\
 & &  \\  
\end{array}
\end{equation}
(The two groups of order $648$ are isomorphic, though differently generated.)

Finally, if $G = [\infty,6,3]$ and $p>3$, we have the same group as in (\ref{grsithth}),   
now yielding a polytope $\cal P$ of type $\{p,6,3\}$.  In particular, if $p=5$, the facets are maps $\{5,6\}_{4}$, with group $S_{5}\times C_{2}$.  When $p=3$, we again obtain the universal polytope (\ref{363sub3univ}), or its dual, depending on choice of diagram.
\medskip

\section{Groups $[k,l,m]$ with a hyperbolic subgroup $[l,m]$}
\label{polfourtwo}

We now consider the crystallographic groups $G = [k,l,m]$, for which $[l,m]$ is of hyperbolic 
type. The groups with a spherical or Euclidean subgroup $[k,l]$ have already occurred in the 
previous section in the dual setting, so we may assume that $[k,l]$ is also of hyperbolic type. 
Then, except possibly for small primes $p$, all three subspaces $V$, $V_{0}$ and $V_{3}$ are non-singular, but $V_{03}$ is non-singular only if $l \neq \infty$. By 
Corollary~\ref{corfourtwo}, if $l = \infty$, then $G^p$ is a C-group 
for $p \geq 3$. On the other hand, by Corollary~\ref{corfourone}, if 
$l \neq \infty$, then $G^p$ can only be a C-group if $p$ is small. 
\smallskip

\noindent
{\bf 5.1.\  The groups $[k,\infty,m]$}

We begin by discussing the groups $G=[k,l,m]$ with $l = \infty$ and with diagrams specified 
below. Inspection of the discriminants shows that $V$ is non-singular, except when 
\[ (k,m,p) \mbox{ or  } (m,k,p) = 
\left\{ \begin{array}{l}
(3,3,7),(4,3,5),(6,3,3),(6,3,13), (4,4,3),(6,4,3),(6,4,7),\\
(6,6,3),(6,6,5),(\infty,6,3).
\end{array} \right.\]
Moreover, $V_{0}$ and $V_3$ are non-singular except occasionally when $p=3$.
If all three spaces $V$, $V_{0}$, $V_{3}$ are non-singular, then the 
vertex-figures of the polytope $\cal P$ with group $G^p$ are the maps of type $\{p,m\}$ or 
$\{p,p\}$ associated with $[\infty,m]^{p}$ (see \cite[5.7,5.9,5.11,5.12]{monsch}), and the facets 
are duals of such maps (with $m$ replaced by $k$). 

When $G=[k,\infty,3]$, with $k=3$, $4$, $6$ or $\infty$, we take the diagram
\[ \stackrel{1}{\bullet}\!\frac{}{\;\;\;\;\;}\!
\stackrel{1}{\bullet}\!\frac{}{\;\;\;\;\;}\!
\stackrel{4}{\bullet}\!\frac{}{\;\;\;\;\;}\!
\stackrel{4}{\bullet}\,,
\quad 
\stackrel{2}{\bullet}\!\frac{}{\;\;\;\;\;}\!
\stackrel{1}{\bullet}\!\frac{}{\;\;\;\;\;}\!
\stackrel{4}{\bullet}\!\frac{}{\;\;\;\;\;}\!
\stackrel{4}{\bullet}\,,
\quad 
\stackrel{3}{\bullet}\!\frac{}{\;\;\;\;\;}\!
\stackrel{1}{\bullet}\!\frac{}{\;\;\;\;\;}\!
\stackrel{4}{\bullet}\!\frac{}{\;\;\;\;\;}\!
\stackrel{4}{\bullet}
\quad \mbox{ or }\quad
\stackrel{4}{\bullet}\!\frac{}{\;\;\;\;\;}\!
\stackrel{1}{\bullet}\!\frac{}{\;\;\;\;\;}\!
\stackrel{4}{\bullet}\!\frac{}{\;\;\;\;\;}\!
\stackrel{4}{\bullet}, \]
respectively. Then $\mbox{disc}(V) \sim -7$, $-5$, $-39$ or $-1$, respectively, so $V$ may be 
singular for $p = 3$, $5$, $7$ or $13$.

From $[3,\infty,3]$ we obtain a self-dual regular $4$-polytope $\cal P$ of 
type $\{3,p,3\}$ for all $p\geq 3$. In the  non-singular cases with $p \neq 3,7$ we have 
\begin{equation}
\label{thinfth}
\Gamma({\cal P}) = 
\left\{ \begin{array}{ll}
O_1(4,p,1)\;,  & \mbox{ if }  p \equiv 1,9,11,15,23,25 \pmod{28}\\
O_1(4,p,-1)\;, & \mbox{ if } p \equiv 3,5,13,17,19,27 \pmod{28}.
\end{array}\right. 
\end{equation}
In particular, if $p=5$, then $\Gamma({\cal P}) = O_1(4,5,-1)$, and $\cal P$ has $130$ icosahedral facets and $130$ dodecahedral vertex-figures. For $p=7$ 
the $7^3$ facets are maps  $\{3,7\}_8$, 
and so dually the vertex-figures are maps $\{7,3\}_8$. Of course, 
for $p=3$ we get the $4$-simplex with group $S_5$.

The group $G = [4,\infty,3]$, with $p > 5$, yields regular polytopes $\cal P$ of type 
$\{4,p,3\}$ with group 
\begin{equation}
\label{foinfth}
\Gamma({\cal P}) = 
\left\{ \begin{array}{ll}
O_1(4,p,1)\;,  & \mbox{ if }  p \equiv 1,7,9,23 \pmod{40}\\
O_1(4,p,-1)\;, & \mbox{ if } p \equiv 17,31,33,39 \pmod{40}\\
O(4,p,1)\;,  & \mbox{ if }  p \equiv 3,21,27,29 \pmod{40}\\
O(4,p,-1)\;, & \mbox{ if } p \equiv 11,13,19,37 \pmod{40} .
\end{array}\right. 
\end{equation}
For $p=3$ we obtain the $4$-cube $\{4,3,3\}$. When $p=7$, the vertex-figures of $\cal P$ are 
isomorphic to Klein's map $\{7,3\}_8$ of genus $3$. Lastly,
when $p=5$ we obtain a polytope with $125$ facets of type $\{4,5\}_6$
and with $250$ dodecahedral vertex-figures.

When $G = [6,\infty,3]$, with $p \neq 3,13$, we obtain regular polytopes $\cal P$ of type 
$\{6,p,3\}$ with group
\begin{equation}
\label{siinfth}
\Gamma({\cal P}) = 
\left\{ \begin{array}{ll}
O_1(4,p,1)\;,  & \mbox{ if }  p \equiv 1,11,25,47,49,59,61,71,83,119,121,133 \pmod{156}\\
O_1(4,p,-1)\;, & \mbox{ if } p \equiv 23,35,37,73,85,95,97,107,109,131,145,155 \pmod{156}\\
O(4,p,1)\;,  & \mbox{ if }  p \equiv 5,41,43,55,79,89,103,125,127,137,139,149 \pmod{156}\\
O(4,p,-1)\;, & \mbox{ otherwise}.\\
\end{array}\right. 
\end{equation}
For $p=5$ we have a polytope with facets isomorphic to $\{6,5\}_{4}$ and with dodecahedral 
vertex-figures. When $p=3$, both the original diagram 
and the alternative

\[ \stackrel{1}{\bullet}\!\frac{}{\;\;\;\;\;\;}\!
\stackrel{3}{\bullet}\!\frac{}{\;\;\;\;\;\;}\!
\stackrel{12}{\bullet}\!\frac{}{\;\;\;\;\;\;}\!
\stackrel{12}{\bullet}\;\;, \]
yield the universal polytope described earlier in (\ref{sithun}). Finally, even though $V$
is singular for $p=13$, we still obtain a C-group $G^{13}$ of order $13^4 \cdot 7 \cdot 4!$.

Suppose now that $G = [\infty,\infty,3]$. 
For  $p \geq 3$, we obtain a regular polytope of type  
$\{p,p,3\}$,  which we shall denote $\mathcal{C}_{p,p,3}$, again
to suggest connections with the classical star-polytopes. To explain this, we resurrect
the earlier group $[\infty, 3,3]$, now with generators  $\langle s_0, s_1, s_2, s_3 \rangle$,
to which we apply the \textit{mixing operation}
$$
(s_0, s_1, s_2, s_3) \rightarrow (s_0, s_1 s_0 s_1, s_2, s_3) =: 
(r_0, r_1, r_2, r_3) \; .
$$
From this construction we find that the reflection group 
$G = [\infty,\infty,3] := \langle r_0, r_1, r_2, r_3 \rangle$ is a subgroup of 
$[\infty, 3,3]$, with index $4$ in characteristic $0$. (See \cite[Lemma 2.1]{sqf}
or \cite{dens} for an equivalent geometric dissection.) 
Now, as explained in the proof of \cite[Thm. 7D16(a)]{arp}, the index collapses to $1$
in characteristic $p \geq 3$. Thus $G^p$ is the group described in
(\ref{caseinf33}), though with the new generators.

Again for  $p=3$, $\mathcal{C}_{3,3,3}$ is the $4$-simplex with group
$S_5$.  If $p=5$, then
\[ \Gamma(\mathcal{C}_{5,5,3}) = O_{1}(4,5,1) \cong H_{4} , \] 
which is isomorphic to the symmetry group of the 
$120$-cell $\{5,3,3\}$ (or any regular star-polytope in 
Euclidean $4$-space associated with it).  
The facets of $\mathcal{C}_{5,5,3}$ are maps $ \{5,5\!\mid\! 3\}$ ( $= \mathcal{M}_{5,5}$) of 
genus $4$, which can be metrically realized in Euclidean $3$-space by the small stellated 
dodecahedron $\{\frac{5}{2},5\}$; the vertex-figures of $\mathcal{C}_{5,5,3}$ are 
dodecahedra $\{5,3\}$ ( $= \mathcal{M}_{5,3}$).  We 
know from \cite[Thm. 7D16]{arp} (or \cite{grstar}) that the universal abstract polytope 
$\{\{5,5\!\mid\! 3\}, \{5,3\}\}$ is isomorphic to the regular star-polytope $\{\frac{5}{2},5,3\}$ 
in Euclidean $4$-space. Therefore, since they have the same group, we must have 
\[ \mathcal{C}_{5,5,3} \cong \{\{5,5\!\mid\! 3\}, \{5,3\}\} \cong {\textstyle \{\frac{5}{2},5,3\}}. \]

More generally,  
our earlier conjecture can be restated as:

\noindent
{\bf Conjecture~2}: For primes $p \geq 3$, $\mathcal{C}_{p,p,3}$ is the universal
polytope of type
$$
\{ \, \mathcal{M}_{p,p}\, , \, \mathcal{M}_{p,3}\, \}\;\;.
$$

Geometrically, it is useful to view this mixing of
new generators for the same group as a \textit{stellation} of the  polygonal faces of
$\mathcal{C}_{p,3,3}$, thereby yielding $\mathcal{C}_{p,p,3}$. In a sense, 
the two polytopes share the same edges, which are,  however, allocated differently 
to form the $p$-gons in $\mathcal{C}_{p,p,3}$. For $p>3$, the two polytopes have
equal numbers of facets, polygons and edges, though rather different numbers of vertices.
\medskip

For the groups $G=[k,\infty,4]$, with $k=3$, $4$, $6$ or $\infty$, our diagrams are
\[ \stackrel{1}{\bullet}\!\frac{}{\;\;\;\;\;}\!
\stackrel{1}{\bullet}\!\frac{}{\;\;\;\;\;}\!
\stackrel{4}{\bullet}\!\frac{}{\;\;\;\;\;}\!
\stackrel{2}{\bullet}\,,
\quad 
\stackrel{2}{\bullet}\!\frac{}{\;\;\;\;\;}\!
\stackrel{1}{\bullet}\!\frac{}{\;\;\;\;\;}\!
\stackrel{4}{\bullet}\!\frac{}{\;\;\;\;\;}\!
\stackrel{2}{\bullet}\,,
\quad 
\stackrel{3}{\bullet}\!\frac{}{\;\;\;\;\;}\!
\stackrel{1}{\bullet}\!\frac{}{\;\;\;\;\;}\!
\stackrel{4}{\bullet}\!\frac{}{\;\;\;\;\;}\!
\stackrel{2}{\bullet}
\quad \mbox{ or }\quad
\stackrel{4}{\bullet}\!\frac{}{\;\;\;\;\;}\!
\stackrel{1}{\bullet}\!\frac{}{\;\;\;\;\;}\!
\stackrel{4}{\bullet}\!\frac{}{\;\;\;\;\;}\!
\stackrel{2}{\bullet}, \]
respectively. Now $\mbox{disc}(V) = -5$, $-3$, $-21$ or $-2$, respectively, so $V$ may be 
singular for the primes $p = 3$, $5$ or $7$.  We  have
already discussed the case $k=3$ in the dual setting.

From $[4,\infty,4]$, with $p>3$, we obtain a self-dual regular $4$-polytope $\cal P$ of type 
$\{4,p,4\}$ with group 
\begin{equation}
\label{foinffo}
\Gamma({\cal P}) = 
\left\{ \begin{array}{ll}
O_1(4,p,1)\;,  & \mbox{ if }  p \equiv 1,7 \pmod{24}\\
O_1(4,p,-1)\;, & \mbox{ if } p \equiv 17,23 \pmod{24}\\
O(4,p,1)\;,  & \mbox{ if }  p \equiv 13,19 \pmod{24}\\
O(4,p,-1)\;, & \mbox{ if } p \equiv 5,11 \pmod{24}.
\end{array}\right. 
\end{equation}
For $p=5$, the polytope $\cal P$ has facets isomorphic to $\{4,5\}_6$ and vertex-figures 
isomorphic to $\{5,4\}_6$. When $p=3$, the space $V$ is singular, and we obtain a regular 
toroid $\{4,3,4\}_{(3,0,0)}$ of rank $4$ with automorphism group 
$G^{3} \cong C_{3}^{3} \rtimes [3,4]$. 

The group $G = [6,\infty,4]$, with $p \neq 3,7$, yields regular polytopes $\cal P$ of type 
$\{6,p,4\}$ with group  
\begin{equation}
\label{siinffo}
\Gamma({\cal P}) = 
\left\{ \begin{array}{ll}
O_1(4,p,1)\;,  & \mbox{ if }  p \equiv 1,23,25,71,95,121 \pmod{168}\\
O_1(4,p,-1)\;, & \mbox{ if } p \equiv 47,73,97,143, 145, 167 \pmod{168}\\ 
O(4,p,1)\;,  & \mbox{ if }  p \equiv 
\left\{ \begin{array}{ll}
5,11,17,19,31,37,41,55,85,89,101,103,\\
107,109,115,125,139,155 \pmod{168}
\end{array}\right.\\
O(4,p,-1)\;, & \mbox{ otherwise.}\\
\end{array}\right. 
\end{equation}
For $p=5$ we have a polytope with facets isomorphic to $\{6,5\}_4$ and 
vertex-figures isomorphic to $\{5,4\}_6$. 
When $p=3$  both the given diagram and its alternate yield groups $G^3$
which are similar in $GL_4(\mathbb{Z}_3)$ to the group of order  $2592$ defined by the diagrams in (\ref{634mod3}). We get the same non-universal
polytope of type $\{ \{6,3\}_{(3,0)}, \{3,4\} \}$.
Lastly, in   the singular case  with $p=7$,  we still obtain a C-group of order 
$2^5 \cdot 3 \cdot  7^4$.

When $G = [\infty,\infty,4]$, with $p>3$, the corresponding polytopes $\cal P$ are of type 
$\{p,p,4\}$ and their groups are once more newly generated versions of those  
described  in (\ref{caseinf34}). For $p=3$ we obtain  the  cross-polytope $\{3,3,4\}$, whose group has index $3$ in $O(4,3,1) \cong F_{4}$. For $p=5$, the polytope $\cal P$ has 
facets isomorphic to $\{5,5\!\mid\! 3\} \cong \{\frac{5}{2},5\}$ 
and vertex-figures isomorphic to $\{5,4\}_6$. 

We now turn to the groups $G=[k,\infty,6]$, with $k=3$, $4$, $6$ or $\infty$.
We have already investigated the cases $k=3$ or $4$ in the dual setting, so we need only consider the diagrams
\[ 
\stackrel{3}{\bullet}\!\frac{}{\;\;\;\;\;}\!
\stackrel{1}{\bullet}\!\frac{}{\;\;\;\;\;}\!
\stackrel{4}{\bullet}\!\frac{}{\;\;\;\;\;}\!
\stackrel{12}{\bullet}
\quad \mbox{or} \quad
\stackrel{4}{\bullet}\!\frac{}{\;\;\;\;\;}\!
\stackrel{1}{\bullet}\!\frac{}{\;\;\;\;\;}\!
\stackrel{4}{\bullet}\!\frac{}{\;\;\;\;\;}\!
\stackrel{12}{\bullet}, \]
for which $\mbox{disc}(V) \sim -15$ or $-3$, respectively.  

From $G = [6,\infty,6]$, with $p \neq 3,5$, we obtain self-dual regular polytopes of type 
$\{6,p,6\}$ with group 
\begin{equation}
\label{sifinfsi}
\Gamma({\cal P}) = 
\left\{ \begin{array}{ll}
O_1(4,p,1)\;,  & \mbox{ if }  p \equiv 1,23,47,49 \pmod{60}\\
O_1(4,p,-1)\;, & \mbox{ if } p \equiv 11,13,37,59 \pmod{60}\\
O(4,p,1)\;,  & \mbox{ if }  p \equiv 17,19,31,53 \pmod{60}\\
O(4,p,-1)\;, & \mbox{ otherwise.}
\end{array}\right. 
\end{equation}
For $p=3$ or $5$, each possible diagram $\Delta(G)$ gives a singular space $V$,
though we still obtain polytopes. 
When $p=5$, $\mathcal{P}$ is self-dual with   facets $\{6,5\}_{4}$ and vertex-figures 
$\{5,6\}_{4}$ and a group of order $2^4 \cdot 3 \cdot 5^4$. For $p=3$, the 
essentially distinct diagrams

\[ \stackrel{4}{\bullet}\!\frac{}{\;\;\;\;\;\;}\!
\stackrel{12}{\bullet}\!\frac{}{\;\;\;\;\;\;}\!
\stackrel{3}{\bullet}\!\frac{}{\;\;\;\;\;\;}\!
\stackrel{1}{\bullet} 
\qquad , \qquad
\stackrel{1}{\bullet}\!\frac{}{\;\;\;\;\;\;}\!
\stackrel{3}{\bullet}\!\frac{}{\;\;\;\;\;\;}\!
\stackrel{12}{\bullet}\!\frac{}{\;\;\;\;\;\;}\!
\stackrel{36}{\bullet} 
\qquad  \mbox{and} \qquad 
\stackrel{3}{\bullet}\!\frac{}{\;\;\;\;\;\;}\!
\stackrel{1}{\bullet}\!\frac{}{\;\;\;\;\;\;}\!
\stackrel{4}{\bullet}\!\frac{}{\;\;\;\;\;\;}\!
\stackrel{12}{\bullet} 
\]
describe exactly the same groups $G^3$ as the diagrams
displayed in (\ref{636othera}) and (\ref{636otherb}). Thus only the first 
and last diagrams yield self-dual polytopes.
\smallskip

When $G=[\infty,\infty,6]$ and $p > 3$, the polytope $\cal P$ is of type $\{p,p,6\}$ and its 
group is given by (\ref{grsithth}).  
For $p=5$, $\cal P$ has facets 
$\{5,5\!\mid\! 3\} \cong \{\frac{5}{2},5\}$ and vertex-figures $\{5,6\}_4$, and 
$\Gamma({\cal P}) = O(4,5,-1)$. For $p=3$ each admissible diagram 
yields the universal polytope 
$\{\{3,3\},\{3,6\}_{(3,0)}\}$ (dual to the polytope described in (\ref{sithun})).

It remains to investigate the group $G=[\infty,\infty,\infty]$, with diagram
\[ \stackrel{4}{\bullet}\!\frac{}{\;\;\;\;\;\;}\!
\stackrel{1}{\bullet}\!\frac{}{\;\;\;\;\;\;}\!
\stackrel{4}{\bullet}\!\frac{}{\;\;\;\;\;\;}\!
\stackrel{1}{\bullet} .\]
Then $\mbox{disc}(V) \sim -1$, independent of $p$. 
In each case we obtain
a regular polytope $\mathcal{C}_{p,p,p}$  of type $\{p,p,p\}$, again
related to a classical star-polytope. Note that 
$\mathcal{C}_{p,p,p}$, along with its   facets and vertex-figures, is self-dual. 
Taking an approach similar to that 
for $\mathcal{C}_{p,p,3}$,
we begin with the group $[\infty, 3,3] = \langle s_0, s_1, s_2, s_3 \rangle$,
to which we apply the mixing operation
$$
(s_0, s_1, s_2, s_3) \rightarrow (s_1, s_0 , s_2 s_1 s_0 s_1 s_2, s_3) =: 
(r_0, r_1, r_2, r_3) \; .
$$
Then in characteristic $0$,
 $G = [\infty,\infty,\infty] := \langle r_0, r_1, r_2, r_3 \rangle$ is a subgroup of 
index $6$ in $[\infty, 3,3]$ (see   \cite[part (6)]{dens}). 
Again following  the proof of \cite[Thm. 7D16(c)]{arp}, we find that 
the index collapses to $1$ in characteristic
$p \geq 3$. Once more $G^p$ is the group described
in (\ref{caseinf33}), though with the new generators.

For $p=3$ we again get the regular simplex $\{ 3,3,3\}$. If $p=5$, 
the facets and vertex-figures are maps 
$\{5,5\!\mid\! 3\}$  ($ = \mathcal{M}_{p,p}$) of genus $4$, 
and \[ \Gamma(\mathcal{C}_{5,5,5}) = O_1(4,5,1) \cong H_{4} .\] 
We know from \cite[Thm. 7D16]{arp} (or \cite{grstar}) that the universal abstract regular 
$4$-polytope $\{\{5,5\!\mid\! 3\},\{5,5\!\mid\! 3\}\}$ is isomorphic to the regular star-polytope 
$\{\frac{5}{2},5,\frac{5}{2}\}$ in Euclidean $4$-space, which also has group $H_{4}$. 
Then it follows that
\[ \mathcal{C}_{5,5,5} \cong \{\{5,5\!\mid\! 3\},\{5,5\!\mid\! 3\}\} \cong 
{\textstyle \{5,\frac{5}{2},5\}}. \] 
This suggests a final variant of our earlier conjecture:

\noindent
{\bf Conjecture~3}: For primes $p \geq 3$, $\mathcal{C}_{p,p,p}$ is the universal
polytope of type
$$
\{ \, \mathcal{M}_{p,p}\, , \, \mathcal{M}_{p,p}\, \}\;\;.
$$
\smallskip

\noindent
{\bf 5.2.\  The groups $[k,l,m]$ with $l=3$, $4$ or $6$}

By Corollary~\ref{corfourone}, the remaining
groups $G = [k,l,m]$ with $l=3$, $4$ or $6$ do not generally 
yield a C-group $G^p$. In particular, if $V$, $V_{0}$, $V_{3}$, $V_{0,3}$ are non-singular and 
$G^p_{0}$, $G^p_{3}$ are of orthogonal type, then $G^p$ fails to be a C-group if $p>13$. Hence only small primes need to be considered.

Every group $G=[k,3,m]$, except $[\infty,3,\infty]$, has a spherical or Euclidean subgroup 
$[k,3]$ or $[3,m]$, so $G$ has already been investigated in the previous section.  

When $G = [\infty,3,\infty]$, with diagram
\[ \stackrel{4}{\bullet}\!\frac{}{\;\;\;\;\;\;}\!
\stackrel{1}{\bullet}\!\frac{}{\;\;\;\;\;\;}\!
\stackrel{1}{\bullet}\!\frac{}{\;\;\;\;\;\;}\!
\stackrel{4}{\bullet} , \]
the three subspaces $V$, $V_{0}$, $V_{3}$ have discriminant $-1$, and 
$\mbox{disc}(V_{0,3}) \sim 3$, so all four subspaces are certainly non-singular for $p>3$.  Moreover, $G^p_{0}$, $G^p_{3}$ are of orthogonal type for $p>3$, and ${\epsilon}(V_{0,3}) = 1$ if and only if $p \equiv 1,7 \bmod 12$.  (Recall that ${\epsilon}(V_{0,3}) = 1$ or $-1$ according as $\mbox{disc}(V_{0,3}) \sim -1$ or $-\gamma$, with $\gamma$ a non-square.)  Now Corollary~\ref{corfourone}(a) shows that $G^p$ can 
only be a C-group if $p \leq 7$. Computations with GAP confirm that $G^p$ is indeed a C-group for $p=3, 5$ or $7$.  The corresponding polytope has type $\{p, 3, p\}$,
and we shall denote it by $\mathcal{C}_{p,3,p}$, again to suggest a connection
with the classical cases.  Even though there are just three polytopes
in this family,  the group $G^p$ is still described by
(\ref{caseinf33}) for any prime $p \geq 3$.

Clearly, when $p=3$ we reacquire the simplex $\{3,3,3\}$.   
For $p=5$ we obtain a self-dual polytope 
$\mathcal{C}_{5,3,5}$ with $120$ dodecahedral facets $\{5,3\}$ and 
$120$ icosahedral vertex-figures $\{3,5\}$, and with 
\[ \Gamma(\mathcal{C}_{5,3,5}) = O_{1}(4,5,1) \cong H_{4} .\] 
We know from \cite[Thm. 7D16]{arp} that the regular star-polytope $\{5,3,\frac{5}{2}\}$ in 
$\mathbb{E}^4$ is isomorphic to the abstract regular polytope $\{5,3,5\!\mid\! 3\}$, which is 
the quotient of the hyperbolic tessellation $\{5,3,5\}$ obtained by imposing the extra relation 
\[ (\rho_{0}\rho_{1}\rho_{2}\rho_{3}\rho_{2}\rho_{1})^{3} = 1 \]
on its group $\langle\rho_{0},\ldots,\rho_{3}\rangle$ (the ``$3$" in the symbol for the quotient 
signifies this relation). It is straightforward to check that the corresponding element 
$r_{0}r_{1}r_{2}r_{3}r_{2}r_{1}$ of $G^5$ indeed has order $3$, so we also have 
\[ \mathcal{C}_{5,3,5} \cong \{5,3,{\textstyle \frac{5}{2}}\} \cong \{5,3,5\!\mid\! 3\} . \]
(In fact, this isomorphism proves directly that $G^5$ is a C-group.)

When $p=7$, the self-dual polytope $\mathcal{C}_{7,3,7}$ has $350$ facets isomorphic to 
$\{7,3\}_{8}$ and $350$ vertex-figures isomorphic to $\{3,7\}_{8}$, and 
$\Gamma(\mathcal{C}_{7,3,7}) = O_{1}(4,7,-1)$, of order $7^{2}(7^{4}-1)$. Thus Klein's dual pair of maps $\{7,3\}_{8}$ and $\{3,7\}_{8}$ of genus $3$ occur as facets and vertex-figures of a self-dual regular $4$-polytope. (In contrast to the case $p=5$, we have no interesting presentation for $G^7$.)

All but three groups $[k,4,m]$ have a spherical or Euclidean subgroup $[k,4]$ or $[4,m]$, 
so their polytopes have already been discussed in the previous section. The three exceptions are $G = [6,4,6]$, $[\infty,4,6]$ (or $[6,4,\infty]$) and $[\infty,4,\infty]$,
with respective diagrams  
\[ \stackrel{3}{\bullet}\!\frac{}{\;\;\;\;\;}\!
\stackrel{1}{\bullet}\!\frac{}{\;\;\;\;\;}\!
\stackrel{2}{\bullet}\!\frac{}{\;\;\;\;\;}\!
\stackrel{6}{\bullet}\,,
\;\quad 
\stackrel{4}{\bullet}\!\frac{}{\;\;\;\;\;}\!
\stackrel{1}{\bullet}\!\frac{}{\;\;\;\;\;}\!
\stackrel{2}{\bullet}\!\frac{}{\;\;\;\;\;}\!
\stackrel{6}{\bullet}
\qquad \mbox{ and }\qquad
\stackrel{4}{\bullet}\!\frac{}{\;\;\;\;\;}\!
\stackrel{1}{\bullet}\!\frac{}{\;\;\;\;\;}\!
\stackrel{2}{\bullet}\!\frac{}{\;\;\;\;\;}\!
\stackrel{8}{\bullet}. \]

When $G = [6,4,6]$, the four subspaces $V$, $V_{0}$, $V_{3}$ and $V_{0,3}$ have discriminants $\sim -7$, $-6$, $-6$ and $1$, respectively, and hence   are non-singular for $p \neq 3, 7$. Moreover, $G^p_{0}$, $G^p_{3}$ are of orthogonal type for $p>3$, and ${\epsilon}(V_{0,3}) = 1$ if and only if $p \equiv 1 \bmod 4$. Then Corollary~\ref{corfourone}(a) implies that $G^p$ can only be a C-group if $p \leq 7$. In fact, $G^7$ also fails to be a C-group.

However, $G^5$ is a C-group by Corollary~\ref{corfourone}(b), 
because $G_{0}^{5}$ and $G_{3}^{5}$ are full orthogonal groups. 
The corresponding regular polytope $\cal P$ is self-dual of type 
$\{6,4,6\}$. A dual pair of 
Petrie-Coxeter polyhedra $\{6,4\!\mid\! 3\}$ and $\{4,6\!\mid\! 3\}$ occur as the types of the 
$130$ facets and $130$ vertex-figures of $\cal P$, respectively; and 
$\Gamma({\cal P}) = O(4,5,-1)$.

When $p=3$, each possible diagram yields a  polytope with
facets isomorphic to $\{6,4\}_{4}$ and vertex-figures 
isomorphic to $\{4,6\}_{4}$. Both the diagram given above for this case
and the alternate

\[ \stackrel{1}{\bullet}\!\frac{}{\;\;\;\;\;}\!
\stackrel{3}{\bullet}\!\frac{}{\;\;\;\;\;}\!
\stackrel{6}{\bullet}\!\frac{}{\;\;\;\;\;}\!
\stackrel{2}{\bullet}\;, \]
yield isomorphic, self-dual polytopes with group order 2592. (This is  
a \textit{new} group of that order, not isomorphic for example to 
the group defined for $p=3$ by the diagrams in (\ref{634mod3}).)
For the remaining diagram 
\[ \stackrel{1}{\bullet}\!\frac{}{\;\;\;\;\;}\!
\stackrel{3}{\bullet}\!\frac{}{\;\;\;\;\;}\!
\stackrel{6}{\bullet}\!\frac{}{\;\;\;\;\;}\!
\stackrel{18}{\bullet}\; , \]
we find that $G^3$ has order $7776$ and  is the automorphism
group of the universal self-dual polytope 
\begin{equation}\label{univ646}
\{ \{6,4\}_{4}\, , \,  \{4,6\}_{4} \}\;\;.
\end{equation}

For $G = [\infty,4,6]$, the discriminants of 
$V$, $V_{0}$, $V_{3}$ and $V_{0,3}$ are   
$\sim -6$, $-3$, $-1$ and $1$, respectively. 
Hence, if $p>3$, then all four subspaces are 
non-singular, $G^p_{0}$, $G^p_{3}$ are of orthogonal type, and ${\epsilon}(V_{0,3}) = 1$ 
if and only if $p \equiv 1 \bmod 4$. Then Corollary~\ref{corfourone}(a) shows again that $G^p$ can only be a C-group if $p \leq 7$. Indeed, for $p=7$ we do obtain a C-group 
of order $2^9 \cdot 3^2 \cdot 7^2$.

We see also that $G^5$ is a C-group by 
Corollary~\ref{corfourone}(b). The 
corresponding polytope $\cal P$
has $120$ facets isomorphic to $\{5,4\}_6$ and $120$ vertex-figures isomorphic to 
$\{4,6\!\mid\! 3\}$, and 
\[ \Gamma({\cal P}) = O(4,5,1) \cong H_{4} \rtimes C_{2} . \]
When $p=3$, both the given diagram and its alternate
produce the dual of the universal polytope
(\ref{univ64sub43}), with group order $2592$.
\smallskip

For $G = [\infty,4,\infty]$, the discriminants of $V$, $V_{0}$, $V_{3}$ and $V_{0,3}$ are
$\sim -2$, $-2$, $-1$ and $1$, respectively, 
so these spaces are non-singular for each $p$. 
Corollary~\ref{corfourone}(a) once again implies that $G^p$ can only be a C-group if 
$p \leq 7$, and  
Corollary~\ref{corfourone}(b) confirms that $G^5$ actually is a C-group.  Moreover, 
computation 
with GAP shows that $G^3$ and $G^7$ also are C-groups. 
The corresponding regular polytopes $\cal P$ are self-dual and are of 
type $\{p,4,p\}$ in each case. 
In particular, for $p=3$ we have 
\[ \Gamma({\cal P}) \cong O(4,3,1) \cong F_{4},\]  
so $\cal P$ is the $24$-cell $\{3,4,3\}$. When $p=5$, the polytope $\cal P$ has $130$ 
facets $\{5,4\}_6$ and $130$ vertex-figures $\{4,5\}_6$, and 
$\Gamma({\cal P}) \cong O(4,5,-1)$. 
Finally, if $p=7$, then $\cal P$ has $350$ vertices and $350$ facets, and 
$\Gamma({\cal P}) = O_{1}(4,7,-1)$; the vertex-figures are isomorphic to the map 
R10.9 of \cite{cond1}.

For the groups $[k,6,m]$ we may assume that $k,m \neq 3$, as the other groups 
have already been
studied  previously.  Now $\mbox{disc}(V_{0,3}) \sim 3$ in each case, so 
${\epsilon}(V_{0,3}) = 1$ if and only if $p \equiv 1,7 \bmod 12$.

When $G = [k,6,4]$, with $k=4$, $6$ or $\infty$, we take the diagram
\begin{equation}
\label{diaksf}
\stackrel{2}{\bullet}\!\frac{}{\;\;\;\;\;}\!
\stackrel{1}{\bullet}\!\frac{}{\;\;\;\;\;}\!
\stackrel{3}{\bullet}\!\frac{}{\;\;\;\;\;}\!
\stackrel{6}{\bullet}\,,
\quad 
\stackrel{3}{\bullet}\!\frac{}{\;\;\;\;\;}\!
\stackrel{1}{\bullet}\!\frac{}{\;\;\;\;\;}\!
\stackrel{3}{\bullet}\!\frac{}{\;\;\;\;\;}\!
\stackrel{6}{\bullet}
\quad \mbox{ or }\quad
\stackrel{4}{\bullet}\!\frac{}{\;\;\;\;\;}\!
\stackrel{1}{\bullet}\!\frac{}{\;\;\;\;\;}\!
\stackrel{3}{\bullet}\!\frac{}{\;\;\;\;\;}\!
\stackrel{6}{\bullet},
\end{equation}
respectively.  Then $\mbox{disc}(V) \sim -2$, $-15$ or $-6$, and 
$\mbox{disc}(V_{3}) \sim -6$, $-2$ 
or $-1$, respectively.  Moreover, $\mbox{disc}(V_{0}) \sim -2$ in each case.  It follows that the 
subspaces $V$, $V_{0}$, $V_{3}$ and $V_{0,3}$ certainly are non-singular for $p>3$, except when $k=6$ and $p=5$; in particular, $G^p_{0}$ and $G^p_{3}$ are then of orthogonal type. Now Corollary~\ref{corfourone}(a) eliminates all primes $p>13$. 

For $G = [4,6,4]$, the subgroups $G_{0}^p$ and $G_{3}^p$ are full orthogonal for all primes 
$p \leq 13$, so Corollary~\ref{corfourone}(b) shows (for $p>3$)
that $G^p$ is   a C-group only when $p=5,7$. For $p=5$, 
we obtain a self-dual regular polytope $\cal P$ with $130$ facets isomorphic to 
$\{4,6\!\mid\! 3\}$ and $130$ vertex-figures isomorphic to $\{6,4\!\mid\! 3\}$, and with group 
$\Gamma({\cal P}) = O(4,5,-1)$. 
Thus $\cal P$ has a dual pair of Petrie-Coxeter polyhedra as its facets and vertex-figures.
For $p=7$, we have another self-dual polytope $\cal P$ of 
type $\{4,6,4\}$, now with $350$  
vertices, $350$ facets, and automorphism group $O(4,7,-1)$.
Finally, when $p=3$, $G^3$ has order $5184$ and is the automorphism group
of the universal, self-dual polytope
\begin{equation}\label{univ464mod3}
\{ \{4,6\}_4 \, ,\, \{6,4\}_4 \}\;\;, 
\end{equation}
with $36$ facets (and vertices).

Let $G = [6,6,4]$. Now Corollary~\ref{corfourone}(b)  applies only 
to $p=7$ among  primes  
$p \leq 13$, so in particular $G^7$ is a C-group. The corresponding regular polytope 
$\cal P$ is of type $\{6,6,4\}$, has $350$ vertices and $350$ facets, 
and has $O(4,7,-1)$ as its 
group. Using GAP we find that $G^p$ is also a C-group for the remaining primes 
$p \leq 13$, thereby giving a few regular polytopes $\cal P$   of type $\{6,6,4\}$. 

For $p=3$, the group $G^3$ depends on the diagram used for   reduction. 
Both the middle diagram in (\ref{diaksf}) and
$$
\stackrel{1}{\bullet}\!\frac{}{\;\;\;\;\;}\!
\stackrel{3}{\bullet}\!\frac{}{\;\;\;\;\;}\!
\stackrel{1}{\bullet}\!\frac{}{\;\;\;\;\;}\!
\stackrel{2}{\bullet}
$$
give the group of order $864$ for the universal polytope 
$$\{ \{6,6\}^\pi_{(1,1)}\, ,\, \{ 6,4\}_4 \}$$
whose $12$ facets are isomorphic to the 
Petrial of Sherk's map $\{6,6\}_{(1,1)}$, and with $6$ vertex-figures 
isomorphic to $\{6,4\}_4$. 
On the other hand, both 
 
\[ \stackrel{1}{\bullet}\!\frac{}{\;\;\;\;\;\;}\!
\stackrel{3}{\bullet}\!\frac{}{\;\;\;\;\;\;}\!
\stackrel{9}{\bullet}\!\frac{}{\;\;\;\;\;\;}\!
\stackrel{18}{\bullet} 
  \quad \mbox{and} \quad
  \stackrel{9}{\bullet}\!\frac{}{\;\;\;\;\;\;}\!
\stackrel{3}{\bullet}\!\frac{}{\;\;\;\;\;\;}\!
\stackrel{1}{\bullet}\!\frac{}{\;\;\;\;\;\;}\!
\stackrel{2}{\bullet} 
\]
yield the group of order $7776$ for the universal polytope
$$\{ \{6,6\}^\pi_{(3,0)}\, ,\, \{ 6,4\}_4 \}\;,$$
whose $36$ facets isomorphic to the Petrial of   $\{6,6\}_{(3,0)}$.
(This is \textit{not} the group for the polytope in (\ref{univ646}).)

When $p=5$, the underlying space $V$ is 
singular and $G^5$ is of order $30000$; the corresponding polytope $\cal P$ has $125$ 
facets isomorphic to the map R11.5 of \cite{cond1}, and $125$ vertex-figures isomorphic to $\{6,4 \!\mid\! 3\}$. Finally, we have the groups $G^{11} = O(4,11,-1)$ and $G^{13} = O(4,13,-1)$, each contributing a regular polytope of type $\{6,6,4\}$, with $2684$ or $4420$ facets, respectively, and with half as many vertices.

Next let $G = [\infty,6,4]$. The primes $5$ and $7$ are easily seen 
to yield C-groups by Corollary~\ref{corfourone}(b); in fact,
$G^p$ is a C-group for the remaining primes $p \leq 13$.
  In particular, when $p=5$, we have a polytope
$\cal P$ with $120$ facets $\{5,6\}_4$ and $120$ vertex-figures $\{6,4\!\mid\! 3\}$, and with 
\[\Gamma({\cal P}) = O(4,5,1) \cong H_{4} \rtimes C_{2} .\]
When $p=7$, $\cal P$ is of type $\{7,6,4\}$ and has $336$ facets and $336$ vertices, and 
$\Gamma({\cal P}) = O(4,7,1)$. 
Similarly, we find that $G^{11} = O(4,11,+1)$ and $G^{13} = O(4,13,-1)$. These
groups yield regular polytopes of type $\{11,6,4\}$ or $\{13,6,4\}$, respectively, with $2640$ or $4420$ facets and half as many vertices.

When $p=3$, the right-most diagram in (\ref{diaksf}) yields the group
$G^3$ of order $432$ for the universal regular polytope 
$$ \{ \{3,6\}_{(1,1)} \, , \, \{ 6,4 \}_4 \} \;\; .$$  
This polytope has with $12$ toroidal facets $\{3,6\}_{(1,1)}$ and $3$ vertex-figures 
$\{6,4\}_{4}$. However,  the alternate diagram 
\[\stackrel{12}{\bullet}\!\frac{}{\;\;\;\;\;\;}\!
\stackrel{3}{\bullet}\!\frac{}{\;\;\;\;\;\;}\!
\stackrel{1}{\bullet}\!\frac{}{\;\;\;\;\;\;}\!
\stackrel{2}{\bullet} \]
gives the group of order $3888$ for the universal regular polytope
$$ \{ \{3,6\}_{(3,0)} \, , \, \{ 6,4 \}_4 \} \;\; ,$$  
now with $36$ toroidal facets $\{3,6\}_{(3,0)}$ and $27$ vertex-figures.

At last then, we are left to consider only  the groups $G = [6,6,6], [\infty,6,6]$ and
$[\infty,6,\infty ]$, with  diagrams

\[ \stackrel{3}{\bullet}\!\frac{}{\;\;\;\;\;\;}\!
\stackrel{1}{\bullet}\!\frac{}{\;\;\;\;\;\;}\!
\stackrel{3}{\bullet}\!\frac{}{\;\;\;\;\;\;}\!
\stackrel{1}{\bullet}\, ,\, 
\stackrel{4}{\bullet}\!\frac{}{\;\;\;\;\;\;}\!
\stackrel{1}{\bullet}\!\frac{}{\;\;\;\;\;\;}\!
\stackrel{3}{\bullet}\!\frac{}{\;\;\;\;\;\;}\!
\stackrel{1}{\bullet} 
  \quad \mbox{and} \quad
  \stackrel{4}{\bullet}\!\frac{}{\;\;\;\;\;\;}\!
\stackrel{1}{\bullet}\!\frac{}{\;\;\;\;\;\;}\!
\stackrel{3}{\bullet}\!\frac{}{\;\;\;\;\;\;}\!
\stackrel{12}{\bullet} \;,
\]
and $\mbox{disc}(V) \sim -11$, $-1$ and $-3$,  respectively. Again
the subspaces $V$, $V_{0}$, $V_{3}$ and $V_{0,3}$ are non-singular for $p>3$
(excluding $p=11$ for the left  diagram); and then
$G^p_{0}$ and $G^p_{3}$ are   of orthogonal type. By Corollary~\ref{corfourone}(a) 
we  exclude all $p>13$. For each diagram, ${\epsilon}(V_{0,3}) = 1$ 
if and only if $p \equiv 1\, , \, 7 \bmod 12$. Thus, by Corollary~\ref{corfourone}(b), 
we do get a C-group for $p = 5, 7$.  In fact, from GAP we find for each diagram 
that $G^p$ is a C-group for $p = 3, 11, 13$, too. We describe some features of the corresponding polytopes in the more notable cases only.

Let us first consider  $G = [6,6,6]$.  The resulting polytope
$\mathcal{P}$ is self-dual for $p = 5, 7, 11, 13$. When $p=5$, the facets
(and vertex-figures) are copies of the map R11.5 listed in \cite{cond1}. 
For $p=11$, $V$ is singular and $G^{11}$ has order $11^4 \cdot 5!$\,. 
For $p=3$, the facets and vertex-figures are the Petrials 
of Sherk's maps, as described earlier. The diagrams 

$$ \stackrel{3}{\bullet}\!\frac{}{\;\;\;\;\;\;}\!
\stackrel{1}{\bullet}\!\frac{}{\;\;\;\;\;\;}\!
\stackrel{3}{\bullet}\!\frac{}{\;\;\;\;\;\;}\!
\stackrel{1}{\bullet} 
  \quad \mbox{ and } \quad
  \stackrel{1}{\bullet}\!\frac{}{\;\;\;\;\;\;}\!
\stackrel{3}{\bullet}\!\frac{}{\;\;\;\;\;\;}\!
\stackrel{9}{\bullet}\!\frac{}{\;\;\;\;\;\;}\!
\stackrel{27}{\bullet}  
$$
yield, respectively, the groups of orders $432$ and $11664$ for the self-dual, universal 
regular polytopes
$$\{ \{6,6\}^\pi_{(1,1)} \, , \, \{6,6\}^\pi_{(1,1)} \}
\quad \mbox{ and } \quad
\{ \{6,6\}^\pi_{(3,0)} \, , \, \{6,6\}^\pi_{(3,0)} \}\;\;.
$$
The remaining pertinent  diagrams
$$ \stackrel{1}{\bullet}\!\frac{}{\;\;\;\;\;\;}\!
\stackrel{3}{\bullet}\!\frac{}{\;\;\;\;\;\;}\!
\stackrel{9}{\bullet}\!\frac{}{\;\;\;\;\;\;}\!
\stackrel{3}{\bullet} 
  \quad \mbox{ and } \quad
  \stackrel{3}{\bullet}\!\frac{}{\;\;\;\;\;\;}\!
\stackrel{1}{\bullet}\!\frac{}{\;\;\;\;\;\;}\!
\stackrel{3}{\bullet}\!\frac{}{\;\;\;\;\;\;}\!
\stackrel{9}{\bullet}  
$$
yield the group of order $1296$ for the universal  regular polytope
$$\{ \{6,6\}^\pi_{(3,0)} \, , \, \{6,6\}^\pi_{(1,1)} \}\;\;,$$
and its dual, respectively.

Next suppose $G = [\infty,6,6]$. For $p=5$ we reobtain 
$\Gamma({\cal P}) = O(4,5,1) \cong H_{4} \rtimes C_{2}$
as the group of order $28800$ for a polytope $\mathcal{P}$
whose $120$ facets are isomorphic to $\{5,6\}_4$ and whose $120$ 
vertex-figures are copies of map R11.5 in \cite{cond1}.

For $p=3$, the situation is quite analogous to that in the
previous case. The diagrams 

$$ \stackrel{4}{\bullet}\!\frac{}{\;\;\;\;\;\;}\!
\stackrel{1}{\bullet}\!\frac{}{\;\;\;\;\;\;}\!
\stackrel{3}{\bullet}\!\frac{}{\;\;\;\;\;\;}\!
\stackrel{1}{\bullet} 
  \quad \mbox{ and } \quad
  \stackrel{36}{\bullet}\!\frac{}{\;\;\;\;\;\;}\!
\stackrel{9}{\bullet}\!\frac{}{\;\;\;\;\;\;}\!
\stackrel{3}{\bullet}\!\frac{}{\;\;\;\;\;\;}\!
\stackrel{1}{\bullet}  
$$
yield, respectively, 
the groups of orders $216$ and $5832$ for the  duals of the first 
and last of the  universal 
regular polytopes displayed in (\ref{univ663}).
The two other pertinent  diagrams
$$ \stackrel{12}{\bullet}\!\frac{}{\;\;\;\;\;\;}\!
\stackrel{3}{\bullet}\!\frac{}{\;\;\;\;\;\;}\!
\stackrel{1}{\bullet}\!\frac{}{\;\;\;\;\;\;}\!
\stackrel{3}{\bullet} 
  \quad \mbox{ and } \quad
\stackrel{4}{\bullet}\!\frac{}{\;\;\;\;\;\;}\!
\stackrel{1}{\bullet}\!\frac{}{\;\;\;\;\;\;}\!
\stackrel{3}{\bullet}\!\frac{}{\;\;\;\;\;\;}\!
\stackrel{9}{\bullet}  
$$
yield isomorphic groups of order $648$. The given generators
then provide the (non-isomorphic) duals of the second and third of the  universal  regular
polytopes  in (\ref{univ663}).

Finally, consider $G = [\infty,6,\infty]$. When $p=3$ we again get the 
universal polytope
described in (\ref{363sub3univ}). For $p = 5, 7,  11, 13$,  we obtain 
self-dual polytopes of type $\{p,6,p\}$,   with groups
$O(4,5,-1)$, $O(4,7,+1)$, $O_1(4,11,-1)$ and $O_1(4,13,+1)$,
respectively. In particular, for $p=5$, we get the universal
self-dual regular polytope 
$$ \{ \{5,6\}_4 \, , \, \{6,5\}_4 \}\;\;.$$

\end{document}